\documentclass[11pt]{article}   
\usepackage{amssymb,amscd,latexsym}   
\usepackage{amsmath}
\usepackage{amsthm}
\usepackage[all]{xy}
\usepackage[dvips]{graphicx}
\usepackage{hyperref}
\newcommand{\rar}{\rightarrow}
\newcommand{\lar}{\longrightarrow}
\newcommand{\llar}{-\kern-5pt-\kern-5pt\longrightarrow}

\newtheorem{Theorem}{Theorem}[section]
\newtheorem{Lemma}[Theorem]{Lemma}

\newtheorem{Proposition}[Theorem]{Proposition}
\newtheorem{Remark}[Theorem]{Remark}
\newtheorem{Example}[Theorem]{Example}

\newtheorem{Definition}[Theorem]{Definition}
\newtheorem{Question}[Theorem]{Question}
\def\sqr#1#2{{\vcenter{\hrule height.#2pt
        \hbox{\vrule width.#2pt height#1pt \kern#1pt
            \vrule width.#2pt}
        \hrule height.#2pt}}}
\def\phi{\varphi}
\def\demo{\noindent{\bf Proof. }}
\def\square{\mathchoice\sqr64\sqr64\sqr{4}3\sqr{3}3}
\def\qed{\hspace*{\fill} $\square$}

\def\XX{{\bf X}}

\def\fm{{\mathfrak m}}

\def\ii{\'{\i}}

\def\hht{{\rm ht}\,}

\def\ker{{\rm ker}\,}

\def\restr{{\kern-1pt\restriction\kern-1pt}}

\def\NN{\mathbb N}

\def\pp{{\mathbb P}}

 2

 3

 4

 5

 4

 3

 2

 1

 1

 2

 3

 1

 0

 1

 2

 3

 4

 5

\begin{document}
\begin{center}
{\Large{\bf\sc Homaloidal nets and ideals of fat points I}}
\footnotetext{Mathematics Subject Classification 2010
 (MSC2010). Primary  13D02, 13H10, 13H15, 14E05, 14E07, 14M05;
Secondary 13A02,  13C14,  14C20.}

\vspace{0.3in}

\hspace{-15pt}{\large\sc Zaqueu Ramos}\footnote{Partially
supported by a CNPq post-doc fellowship (151229/2014-7).}
\quad\quad\quad
 {\large\sc Aron  Simis}\footnote{Partially
supported by a CNPq grant (302298/2014-2) and by a CAPES-PVNS Fellowship.}

\end{center}


\bigskip

\begin{abstract}

One considers plane Cremona maps with proper base points and the {\em base ideal} generated by the linear system of forms defining the map.
The object of this work is  the interweave between the algebraic  properties of the  base ideal  and those of the ideal of these points fattened by the virtual multiplicities arising from the linear system.
One reveals  conditions which naturally regulate this association, with particular emphasis on the homological side.
While most classical numerical inequalities concern the three highest virtual multiplicities, here one emphasizes also the role of one single highest multiplicity.
In this vein one describes classes of Cremona maps for large and small value of the highest virtual multiplicity.
One also deals with the delicate property as to when the base ideal is non-saturated and the structure of its saturation.

\end{abstract}

\section*{Introduction}

Let $k$ be an infinite field and let
$R=k[\XX]=k[X_0,\ldots,X_n]$ denote a polynomial ring over $k$,
endowed with the ordinary standard grading.
A rational map $F\colon\pp^n\dasharrow \pp^n$ is defined by a linear system of
$n+1$ independent forms of the same degree  in $R$.
If $F$ is birational then it is  called a Cremona map.
Beyond the notable modern development on the structure of the Cremona group, a body of results on the nature, structure and uses of an individual Cremona map has lately come up that draws on modern geometric and algebraic  tools (for a very short sample, see \cite{Mella}, \cite{Blanc}, \cite{NewCremona}, \cite{CRS}, \cite{Ciliberto-et-al}, \cite{Dumnicki}, \cite{Pan}, \cite{PanStellar}, \cite{PanRusso},
\cite{bir2003},  \cite{SimisVilla}, \cite{{CremonaMexico}}).
This paper goes along this line and focus only on the plane case (i.e., $n=2$).

The core of the work is a study of the intertwining between the algebraic properties of the base locus of $F$, here taken in the sense of the strict {\em base ideal}
$I\subset R$ generated by the linear system of forms defining $F$, and those of the ideal $I(\mathbf{p},\boldsymbol{\mu})$ of general fat points $\mathbf{p}$ with the virtual multiplicities $\boldsymbol{\mu}$ of this linear system.
A first question is how homaloidal virtual multiplicities imply special properties of the corresponding ideal $J\subset R$ of fat points in terms of its Hilbert function and minimal free resolution.
Conversely, it is natural to ask how to recover information about the base ideal of the Cremona map out of the features of $J$, such as, for example, the homological structure of the $1$-dimensional ring $R/I$.

When relating Cremona maps (or, more generally rational maps of $\pp^2$) to ideals of fat points, a technical difficulty arises at the outset, and that is the fact that the defining linear system of a Cremona map determines along with base points in $\pp^2$ also infinitely near base points.
Capturing the entire picture in the passage to the current notion of fat points becomes flawed.
To partially overcome this problem, it has been proposed in \cite{HS} to introduce a new ideal, to wit, given a Cremona map of degree  $d\geq 1$, with {\em cluster} $\mathcal{K}$ of points and virtual multiplicities, one lets $\ell_{\mathcal{K}}(d)\subset R_d$ consist of all forms $f$ of degree $d$ such the
curve $C=V(f)$ passes virtually through $\mathcal{K}$, and  sets
$I_{\mathcal{K}}:=\oplus_{d\geq 0} \ell_{\mathcal{K}}(d)\subset R$, a nice homogeneous ideal.
As it turns out, there are inclusions $I\subset I_{\mathcal{K}}\subset 
I(\mathbf{p},\boldsymbol{\mu})$, where $\mathbf{p}$ are the proper points of the map and $\boldsymbol{\mu}$ are the corresponding virtual multiplicities. Moreover, the rightmost inclusion is an equality if the map has only proper points (i.e., a {\em simple} Cremona map in the classical terminology \cite{alberich}). Of course, the hardship to capture this construction in purely algebraic terms is the 
notion of ``passing virtually through $\mathcal{K}$'' -- for the details see  \cite[Introduction]{HS}.

A similar notion has been given by Harbourne in \cite{Harb1}
in terms of the push-down of a suitably defined divisor on the blowup of points on $\pp^2$.
If there is an explicit ideal of $R$ to be extracted from the latter notion it ought to be very similar, if not the same, as
the one in \cite{HS}.

Dealing with these additional data would certainly  cloud the central idea of this paper. Therefore, we assume that base points are points of $\pp^2$ and, as a rule, general points.
Since Cremona maps densely exist in this context of points, for any proper homaloidal type, we will be dealing with a vast universe.
In other words, we will emphasize the role of the proper homaloidal type and of a Cremona map with this type and based on proper points only.
Of course, any set of points with random coordinates is considered to be general. In the present context the choice of sufficiently general points could be more precisely described in terms of avoiding the zero set of minors of a certain matrix -- see, e.g., \cite{Rick}. For ordinary computer calculation, by and large random  point coordinates are taken, though special scripts are available for the second alternative (see, e.g., \cite{HarbGuar}).

For plane Cremona maps of degree at least $2$, asking when the base ideal $I\subset R$ is saturated is tantamount to asking when
$R/I$ is a Cohen--Macaulay ring, i.e., when $I$ is generated by the
maximal minors of a $3\times 2$ homogeneous matrix with entries in $R$.
We will often imprecisely refer to a {\em saturated} or {\em Cohen--Macaulay} Cremona map in such a situation.
While there are many rational maps on $\pp^2$, in any degree, defined by such ideals, the question for Cremona
maps becomes tighter.

As is usually taken for granted while discussing this sort of questions, we assume throughout
that the base field $k$ is algebraically closed (a few times of characteristic zero.)

\smallskip

We now proceed to a more detailed description of the sections.

Section 1 is a brief account of the preliminaries and the terminology; the subsection titles speak for themselves. We have tried to be as self-contained as possible in terms of the usage of these two aspects.

\smallskip

Section 2 contains the main theorem that bridges up between the base ideal $I\subset R$ of a plane  Cremona map of degree $d\geq 4$ with proper points and the associated ideal $J\subset R$  of fat points. The minimal free $R$-resolution of $J$ is made explicit and it is shown that the linear system $J_{d+1}$ defines a birational map onto its image.
The proofs of both results are entirely ideal theoretic, but we have added an explanation, through the calculation of the regularity and a typical property of Cremona maps, as to how the second of these results is included in a theorem of Geramita-Gimigliano-Pitteloud in \cite{GGP}.

We also present a substantial list of homological results relating the two ideals $I\subset J$. For example, it is explained how to get to the free resolution of $I$ from the free resolution of $J$.
In addition, it is shown that if the $R$-module $J/I$ has dimension $1$ then the Cohen--Macaulayness of $R/I$ is tantamount to that of $J/I$, which is a simplification in the theory since the latter is a linearly presented $R$-module.
On the other hand, it is shown that the condition $\dim (J/I)=0$ is equivalent to $J$ being the saturation of $I$, which in turn, in degree $\geq 4$, can only happen for the homaloidal types $(4;2^3,1^3)$ and $(5;2^6)$.

\smallskip

In Section 3 we start a classification of homaloidal types according to the highest virtual multiplicity $\mu=\mu_1$.
The permeating idea is the expectation that for $\mu_1\geq d-4$ any Cremona map on general points of such a homaloidal type is saturated.
While this is sufficiently known for $\mu=d-1$ (de Jonqui\`eres case) we prove this for $\mu= d-2$  using the apparatus of the previous section. This is obtained by showing that the linear matrix $\mathcal{L}$ presenting the ideal $J/I$ has Fitting $I_1(\mathcal{L})$ of minimal codimension.
Alas, we also show, to our regret, that this condition fails for $\mu\leq d-3$.

There are however a set of evidences pointing out to saturation in the cases where $d-4\leq \mu\leq d-3$.
As a matter of fact, we are naturally led to (loosely) conjecture that the Cremona map on general (proper) points is saturated unless $\mu\leq \left \lfloor{d/2}\right \rfloor$.

The last section is a detailed study of a strong positive case of the above conjecture, namely, the case of a homaloidal type $T$ whose virtual multiplicities are all even.
This consideration unexpectedly stretches to a close look at the type $(1/2)T$ obtained from $T$ by dividing all multiplicities by $2$.
We convey that although $(1/2)T$ is not a homaloidal type it nevertheless satisfies similar ``equations of condition'' -- the terminology {\em sub-homaloidal} has been introduced to designate a set of virtual multiplicities satisfying these new numerical conditions.

One notes that passing to the corresponding ideals of fat points for $T$ and  $(1/2)T$, the first of these ideals is the second symbolic power of the second ideal. This in itself allows for a crosscut between properties of the two and sets up a relation between the respective initial degrees.
A finer analysis  will be developed in a sequel to this paper,  where one subverts the order of priority, enhancing properties of the second of these ideals.

This part culminates with a proof that the base ideal $I$ of a Cremona map $\mathfrak{F}:\pp^2\dasharrow\pp^2$ of degree $d$, whose virtual multiplicities are all even, is not saturated under the hypothesis that there are three highest virtual multiplicities equal to $\left \lfloor{d/2}\right \rfloor$.
In addition, one gives a minimal set of explicit generators concocted from the respective equations of three irreducible fundamental curves relative to base points of $\mathfrak{F}^{-1}$.
Finally, one proves that the saturation of $I$ requires only one additional generator and this generator can be made as explicit as possible within the context.

\smallskip

As a pointer, the following are the main results of the paper: Theorem~\ref{MAIN}, Proposition~\ref{Supplement-MAIN},  Proposition~\ref{highest_is_2less}, Proposition~\ref{highest_is_3less}, Proposition~\ref{highest_is_4less}, Proposition~\ref{initial_degree_of_divided}, Theorem~\ref{Main-non-saturated} and Proposition~\ref{resolving_non_saturated}.

\section{Recap of terminology}\label{recap}

In this section we briefly wrap-up the basics of linear systems and of ideals of fat points that are going to be used in the sequel.
For complete details, see \cite{alberich} or  \cite[Introduction]{HS}.

\subsection{Linear systems, base points and Cremona maps}

By a {\em linear system} of plane curves of degree $d$ we mean a $k$-vector subspace $L_d$ of the vector space of forms of degree $d$ in the standard graded polynomial ring 
$R: =k[x,y,z]$. 
A linear system $L_d$ defines a rational map $\mathfrak{L}_d$ with source $\pp^2$ and target $\pp^r$, where $r+1$ is the vector space dimension of $L_d$.
As such, both $L_d$ and its closely associated $\mathfrak{L}_d$ have a fixed
part and a set of base points.
More algebraically, one sets $I:=(L_d)\subset R$,  the ideal generated by the elements of $L_d$.
Write $I=(I:F)\,F$, where $F\in R$ is up to nonzero scalars a uniquely defined form of degree $\leq d$ such that $I:F$ has codimension $\geq 2$.
In other words, $F$ is the $\gcd$ of a set of generators of $I$ and $I:F=(1/F)\,I$.
Then $F$ (or its zero set in $\pp^2$) is the {\em fixed part} of the system, while $V(I:F)\subset \pp^2$ is a finite set since $I:F$ is an ideal of codimension $\geq 2$ -- the elements of $V(I:F)$ are the {\em proper base points} of $L_d$ or of the corresponding rational map.
We say that the linear system has no fixed part or is without fixed part meaning that $\deg(F)=0$ or, equivalently, that $I$ codimension $\geq 2$.

For the purpose of this work,  throughout $I:F$ has codimension exactly $2$ as otherwise the set of base points is empty.
Then $I:F$ has a primary decomposition whose components are
primary ideals associated to the minimal primes $\{P_1,\ldots,P_n\}$ of $R/I:F$ defining the base points
$V(I:F)=\{p_1,\ldots,p_n\}$.
Since $k$ is algebraically closed, every one of these primes is generated by two independent linear forms.
We remark that in the tradition of rational maps the fixed part is often set aside because the generators of $I$ essentially define the same map as the generators of $I:F$. However, in the theory of linear systems it is often the case that linear systems without fixed part have subsystems with nontrivial fixed part. 

For the sake of self-containment we recall the notion of multiplicity or local vanishing order.
Namely, given a variety $X$, a smooth point $p\in X$ and a hypersurface (divisor) $D$ then the (effective) {\em multiplicity} 
$e_p(D)$ of $D$ at $p$
is the order of vanishing of a local equation of $D$ at $p\,$; algebraically, if $f$ is a local equation of $D$ at $p$,
then 
$$e_p(D)=\max\{s\geq 0\,|\, f\in \mathfrak{m}_p^s\},$$ 
where $\mathfrak{m}_p$ is the maximal ideal of the local ring
of $X$ at $p$.
For our purpose, $X$ will always be a projectively embedded surface and, in fact, mostly  $\pp^2$ or the blowup of a set of points therein.

The  {\em virtual multiplicity} of $L_d$ at one of its proper base points $p$ is the nonnegative integer
$$\mu_{p}=\mu_{p}(L_d):=\min\{ e_{p}(f)\,|\,f\in L_d\}=
\min_{f\in L_d}\biggl\{s_f\geq 0\,\big |\, f\in ({P}_{P})^{s_f}\setminus ({P}_{P})^{s_f+1}\biggr\}.$$
(Note that requiring the global containment $f\in {P}^{s_f}\setminus {P}^{s_f+1}$ 
in the definition is equivalent since each $P$ is a complete intersection).
Clearly, the effective multiplicity of an element of $L_d$ at a proper base point $p$ is at least $\mu_{p}$,
while the subset of $L_d$ whose elements have  effective multiplicity equal
to $\mu_{p_j}$ forms an open set $U_{p_j}$ (in the set of coefficients as parameters).
Also, $L_d$ admits infinitely near base points, in contrast to its proper base points, and the virtual multiplicity at such points can be defined in a similar way by blowing up proper points.

\medskip

An important property of a plane Cremona map is that its virtual multiplicities satisfy the classical  {\em equations of condition}  (see \cite[2.5]{alberich}):
\begin{equation}\label{eqs_condition}
\sum_{p}\mu_p=3d-3,\; \sum_{p}\mu_p^2=d^2-1,
\end{equation}
where $p$ runs through the set of (proper and infinitely near) base points of the corresponding $L_d$ with respective  multiplicities $\mu_p$.
An abstract configuration $(d\,;\mu_1,\ldots,\mu_r)$ satisfying the equations of condition is called a {\em homaloidal type}.
A homaloidal type is called {\em proper} if there exists a plane Cremona map with this type. There is an important practical  tool to test whether a given homaloidal
type is proper - it is called {\em Hudson test} (\cite[Corollary 5.3.2]{alberich}).
This idea has  been largely developed by M. Nagata in  \cite{Nagata}.
Geometrically, the test can simply be translated into  composing the rational map with the standard quadratic map based on three points of highest virtual multiplicities. Yet, the encoded arithmetic and algebra keep a few surprises.

\subsection{Ideals of fat points}

Besides the interest in the role of linear systems in the study of rational maps, there is a vast literature on the largest possible (complete) such systems for a fixed degree, the so-called 
 {\em ideals of fat points}.

For convenience, recall the definition: one is given a set $\mathbf{p}=\{p_1,\ldots,p_n\}\subset \pp^2$ of points and respective appended multiplicities
$\boldsymbol\mu=\{\mu_1,\ldots,\mu_n\}$, $\mu_j\in \NN_+$.
Let $L_d(\mathbf{p};\boldsymbol\mu)$ (for a lack of shorter notation) denote the set of all $d$-forms $f\in R$ such that
$e_{p_j}(f)\geq \mu_j$ o for every $j$.
Because $e_{p_j}(f)$ has the properties of a valuation, it follows immediately that this set is a vector
space, i.e., it is a linear system of plane curves of degree $d$.
By definition, it is the largest possible such a linear system with $\mathbf{p}$ as its set of proper base
points and $\mu_j$, $1\leq j\leq n$ as the corresponding virtual multiplicities, and actually it contains
as a vector space all such vector subspaces.
If we let $d$ run through $\NN$, we obtain the $\boldsymbol\mu$-fat ideal of $\mathbf{p}$:
\begin{equation}\label{big_fat}
I(\mathbf{p};\boldsymbol\mu):=\bigoplus_{d\in \NN} L_d(\mathbf{p};\boldsymbol\mu)\simeq \bigcap_{j=1}^n P_j^{\mu_j},
\end{equation}
where $P_j\subset R$ is the homogeneous prime ideal of the point $p_j$.

Clearly, the $\boldsymbol\mu$-fat ideal has codimension $2$.
However, the ideal  $(L_d(\mathbf{p};\boldsymbol\mu))\subset R$ generated in a  single degree $d$ is often of codimension $1$.
If the points $\mathbf{p}$ are  general then some authors write $I(n;\boldsymbol\mu)$ instead, a notation we avoid
due to  the diversity of
 notions of {\em general position}.
We will say {\em general points} (or even {\em sufficiently general points}) in the sense that 
they avoid a proper closed subset of $\pp^2\times\cdots\times \pp^2$ ($n$ factors). 
This vague notion can be made more precise (see, e.g., \cite{Rick}). In any case, for the explicit computation one
can always take a random set of points.

\medskip

The ideal $I(\mathbf{p};\boldsymbol\mu)$ is supple in some senses.
Firstly, it is a perfect ideal, i.e., the ring $R/I(\mathbf{p};\boldsymbol\mu)$ is Cohen--Macaulay.
Secondly, the degree of the corresponding scheme is automatic and depends only on the number of
points and their appended multiplicities:
\begin{equation}\label{multplicity_of_fat}
e(R/I(\mathbf{p};\boldsymbol\mu))=\sum_{j=1}^n\, \frac{\mu_j(\mu_j+1)}{2}.
\end{equation}
Thirdly, it is an integrally closed ideal.
Further, quite a bit is known about the Hilbert function of $R/I(\mathbf{p};\boldsymbol\mu)$ and its
graded Betti numbers.
More generally, we have the following.

Let $J\subsetneq R:=k[x,y,z]$ denote an unmixed homogeneous ideal of codimension $2$.
Since $R/J$ is Cohen--Macaulay the minimal graded $R$-resolution of $J$ has the form
\begin{equation}\label{resolution}
0\rar F_1\stackrel{\phi}{\lar} F_0\lar J\rar 0,
\end{equation}
where $\phi$ is a homogeneous map and $F_0,F_1$ are free graded modules.

Let $h_{R/J}(t):=\dim_k(R_t/J_t)\, (t\geq 0)$ stand for the Hilbert function of $R/J$. Recall that the Hilbert polynomial of $R/J$ is of degree $0$ and coincides with the multiplicity $e(R/J)$.

The following facts are well-known:
\begin{enumerate}
\item[{\rm (A)}] $h_{R/J}(t)$ is strictly increasing till it reaches its maximum $e(R/J)$ and thereon stabilizes.
\item[{\rm (B)}] The least $t$ such that $h_{R/J}(t)=e(R/J)$ is the {\em regularity index} of $R/J$
\item[{\rm (C)}] Set, moreover, $F_0=\oplus_{t>0} R(-t)^{n_t}$, so that the minimal number of generators of $J$ is $\sum_{t>0} n_t$.
Then $n_t=\dim_k{\rm coker}(R_1\otimes_k J_{t-1}\lar J_{t})=\dim_k(J_{t})-\dim_k (R_1J_{t-1})$, where $R_1J_{t-1}$ is short for the vector $k$-subspace of $J_{t}$ spanned by the products of $x,y,z$ by the elements of a basis of $J_{t-1}$ and $R_1\otimes_k J_{t-1}\lar J_{t}$ denote the natural multiplication map.
\end{enumerate}

\section{Ideals of fat points coming from homaloidal types}

Recall that a net of degree $d$ is a linear system $L_d$ spanned by $3$ independent forms in $R_d=k[x,y,z]_d$ without proper common factor.
One says that the net defined by $L_d$ is {\em complete} if $L_d=J_d$, where $J=I(\mathbf{p};\boldsymbol\mu)$ is the ideal of fat points based on the cluster $(\mathbf{p};\boldsymbol\mu)$ of proper points defined by $L_d$.
Note that completeness has to do with a fixed degree $d$, for which it implies equalities $(L_d)=(I_{\mathcal{K}})_d= 
I(\mathbf{p},\boldsymbol{\mu})_d$, in the notation explained in the Introduction.
One says that the net $L_d$ is {\em homaloidal} if it defines a Cremona map of $\pp^2$.

The moral of this section is supported by the following basic result:

\begin{Proposition}\label{completeness}
{\rm (char$(k)=0$)}
Let $L_d$ denote a homaloidal net and $T: 
=(d;\mu_1,\ldots,\mu_r)$ the corresponding proper homaloidal type.
Then:
\begin{enumerate}
\item[{\rm (a)}] $L_d$ is complete.
\item[{\rm (b)}] There exists
a non-empty Zariski open subset  $U\subset \underbrace{\pp^2\times\cdots\times\pp^2}_{r}$
such the coordinates of any tuple of points in $U$ is the set of the base points of a simple Cremona map with  homaloidal type $T$ -- i.e., informally speaking, choosing any set of $r$ ``general'' points in $\pp^2$ will give rise to a complete net as in {\rm (a)}.
\end{enumerate}
\end{Proposition}
\demo
(a) See \cite[Proposition 2.5.2]{alberich}. A more conceptual proof follows by the Enriques criterion  \cite[Proposition 5.1.1]{alberich}, as follows: let $L_d$ span a homaloidal net of degree $d$ defining the given Cremona map.
Set $J:=I(\mathbf{p};\boldsymbol\mu)$.
Clearly, $L_d\subset J_d$. Since $L_d$ is the linear system of a 
Cremona map, its generic member $f$ is irreducible. But $f\in J_d$, hence the Enriques criterion implies that $J_d$ is a homaloidal system, i.e., the linear system of a Cremona map; in particular, $\dim_k (J_d)=3$, so we are done.

(b) See \cite[Theorem 5.2.19]{alberich}).
\qed

\begin{Remark}\rm
One cannot possibly be fooled by the apparent amplitude of Enriques criterion: the existence of irreducible curves in the linear system is an expression of Bertini's theorem, even without the assumption that the type $T=(d;\boldsymbol\mu)$ be homaloidal -- see  \cite[Proposition 9]{GerOre}) for cases in which in fact every member is irreducible.
\end{Remark}

\subsection{Main theorem}

In this part we discuss the homological nature of the ideal of fat points of homaloidal type.

We start with a basic result.

\begin{Lemma}\label{initial_degree}
Let $\mathfrak{F}:\pp^2\dasharrow\pp^2$ be a Cremona map with homaloidal type $(d;\boldsymbol\mu)$ whose base points $\mathbf{p}$ are proper {\rm (}i.e, $\mathbf{p}\subset \pp^2${\rm )} and let $J:=I({\bf p},\boldsymbol\mu)\subset R:=k[x,y,z]$ denote the associated ideal of fat points.
Then:
\begin{enumerate}
\item[{\rm (i)}] The initial degree of $J$ is $d$.
\item[{\rm (ii)}] The Hilbert function of $R/J$ is maximal with $h_{R/J}(t)=e(R/J)$ for $t\geq d$ {\rm (}i.e., the regularity index of $R/J$ is $d${\rm )}
\item[{\rm (iii)}] $n_t=0$ for $t\geq d+2$.
\end{enumerate} 
\end{Lemma}
\demo (i) For this assertion, suppose that $0\neq f\in J_{d-1}$. Then $\{xf,yf,zf\}$ is a $k$-linearly independent subset of $J_d$. Since $\mathfrak{F}$ is a Cremona map, its underlying $2$-dimensional linear system $L_d$ is complete by Proposition~\ref{completeness}, i.e., $L_d=J_d$. Therefore, $\{xf,yf,zf\}$ is a $k$-basis of $L_d$, contradicting the fact that $L_d$ has no fixed part.

(ii) For the second assertion we recall that the equations of conditions satisfied by a Cremona map plus  knowing that $\dim_k(J_d)=\dim_k(L_d)=3$ imply the equality
\begin{equation}\label{eqs_of_condition}
\dim_k(J_d)=\frac{(d+2)(d+1)}{2}-\sum \frac{(\mu_i+1)\mu_i}{2}.
\end{equation}
Since $d$ is the initial degree of $J$, the Hilbert function of $R/J$ is maximal.

(iii) This follows from (ii) (see \cite{DGM}).
\qed

\medskip

We will need the following structural result which could not be found in this exact form in the literature:

\begin{Proposition}\label{one-linear}
Let $\mathfrak{F}:\pp^2\dasharrow\pp^2$ be a Cremona map whose base ideal is perfect and has  a  linear syzygy. Then $\mathfrak{F}$ is a de Jonqui\`eres map.
\end{Proposition}
\demo
Since any quadratic Cremona map is de Jonqui\`eres, we may assume at the outset that $\mathfrak{F}$ has degree $\geq 3$.
Consider the presentation matrix of its base ideal $I$:
$$L:=\begin{pmatrix}
\ell_1 & q_1 \\   
\ell_2 & q_2 \\
\ell_3  &  q 
\end{pmatrix},
$$
where $\ell_i$ is a linear form and $q_i,q$ are forms of degree $\geq 2$.
By \cite[Proposition 3.4]{AHA}, $I$ is not an ideal of linear type. It follows that the codimension of the ideal generated by the entries of $L$ is $\leq 2$ since otherwise $I$ would be generically a complete intersection; but an almost complete intersection which is generically a complete intersection must be of linear type (\cite[Proposition 3.7]{SV1}). In particular,  the three linear forms along the first column cannot be $k$-linearly independent. By suitable elementary operation, we may assume that $\ell_1=x,\ell_2=y, \ell_3=0$.
Then one has $I=(xq,yq,\ell_1q_2-\ell_2q_1)$.
Moreover, by \cite[Corollaire 2.3]{PanStellar} $q$ and $\ell_1q_2-\ell_2q_1$ are relatively prime $z$-monoids.
This means that  $\mathfrak{F}$ is a de Jonqui\`eres map (see \cite[Introduction]{PanSi}).
\qed

\begin{Remark}\rm
There are examples of non Cremona dominant rational maps  $\pp^2\dasharrow\pp^2$ whose base ideal $I$ is perfect and has a unique linear syzygy $($up to scalars$)$. Such examples are easily found where $I$ an ideal of linear type -- equivalently, such that the three coordinates of the linear syzygy are $k$-linearly independent.
\end{Remark}

\smallskip

We are now ready for the main result of this part.

\begin{Theorem}\label{MAIN}
Let $\mathfrak{F}:\pp^2\dasharrow\pp^2$ be a Cremona map with homaloidal type $(d;\boldsymbol\mu)$, with $d\geq 4$, whose base points $\mathbf{p}$ are proper {\rm (}i.e, $\mathbf{p}\subset \pp^2${\rm )} and let $I\subset R$ denote its base ideal. Let $J:=I({\bf p},\boldsymbol\mu)$
denote the associated ideal of fat points.
Then:
\begin{enumerate}
\item[{\rm (a)}] The  minimal graded free resolution of $J$ is of the form
\begin{equation}\label{res_DJ}
0\rar R(-(d+2))^{d-2}\oplus R(-(d+1))^s \stackrel{\varphi}\rar R(-(d+1))^{d-4+s}\oplus R(-d)^3\rar J\rar 0
\end{equation}
where $s$ is the number of independent linear syzygies of the ideal $I$. 
Furthermore, $0\leq s\leq 1$, while $s=1$ if and only if  $\mathfrak{F}$ is a de Jonqui\`eres map.

\item[{\rm (b)}] The linear system $J_{d+1}$ defines a birational mapping of $\pp^2$ onto the image in $\pp^{d+4+s}$.
\end{enumerate}
\end{Theorem}
\demo  (a) By Lemma~\ref{initial_degree} (iii), we know that $J$ has only minimal generators in degrees $d$ and $d+1$ and that in degree $d$ the minimal number of generators is $3$.
Thus, we only have to determine the minimal number of generators $n_{d+1}$.

By Lemma~\ref{initial_degree} (ii) and by the equations of condition one has
$$\dim_k(J_{d+1})=\frac{(d+3)(d+2)}{2}-\sum \frac{(\mu_i+1)\mu_i}{2}
= d+5.
$$
Let $N:=\ker(R_1\otimes_k J_d\lar J_{d+1})$ and $K:={\rm coker}(R_1\otimes_k J_d\lar J_{d+1})$.
From the exact sequence
$$0\rar N\lar R_1\otimes_k J_d\lar J_{d+1}\lar K\rar 0,
$$
and since $\dim_k(R_1\otimes_k J_d)=9$,
we then get $n_{d+1}=d-4+\dim_k(N)$.

We claim that $\dim_k(N)$ coincides with the number of independent linear syzygies of the ideal $I$.
For this to be seen, let $\{D_1,D_2,D_3\}$ denote a $k$-vector basis of $J_d$. Letting $\{T_1,T_2,T_3\}$ denote indeterminates over $R$, one has $k[T_1,T_2,T_3]_1=kT_1+kT_2+kT_3\simeq J_d$ as $k$-vector spaces.
Therefore, $R_1\otimes_k J_d\simeq R_1\otimes_k k[T_1,T_2,T_3]_1\simeq R[T_1,T_2,T_3]_{(1,1)}$ as $k$-vector spaces, where the latter is the graded part in bidegree $(1,1)$ of the standard bigraded polynomial ring $R[T_1,T_2,T_3]$.
With this identification it follows immediately that $N$ coincides with the bigraded $(1,1)$-part of the ideal generated by the bigraded relations of $\{D_1,D_2,D_3\}$ -- the so-called Rees ideal of $I=(D_1,D_2,D_3)$ -- which is identified with the $k$-vector space spanned by the $R$-linear syzygies of $I$. 

This shows the first statement in (a).

For the supplementary statement, let $s$ denote as in the statement  the number of independent linear syzygies of $I$.
Clearly, this is exactly the number of independent $k$-linear relations among the elements $\{xD_j,yD_j,zD_j\}$ spanning $R_1J_d$.
Thus, to show that $0\leq s\leq 1$ it suffices to prove that these forms span a $k$-vector space of dimension at least $8$.
For this purpose, consider a linear form $\ell\in R$ which is a nonzerodivisor on $R/J$. We claim that $xD_1,\,xD_2,\, xD_3,\,\ell D_1,\;\ell D_2,\;\ell D_3$ are $k$-linearly independent.
Indeed, from a relation
 $$\alpha_1xD_1+\alpha_2xD_2+\alpha_3 xD_3+\beta_1\ell D_1+\beta_2\ell D_2+\beta_3\ell D_3=0,$$
 with $\alpha_j,\beta_j\in k$,
one deduces
$$x(\alpha_1D_1+\alpha_2D_2+\alpha_3 D_3)=\ell(-\beta_1D_1-\beta_2D_2-\beta_3 D_3);$$
and hence 
$$(-\beta_1D_1-\beta_2D_2-\beta_3 D_3)/x=f\in k[x,y,z],$$  with $f$ a form of degree $d-1$. Noting that $\ell f\in (D_1,D_2,D_3)\subset J$, since $\ell$ is regular modulo $J$ we must have $f\in J.$ But the initial degree of $J$ being $d$, we conclude that $f=0.$ 
Bringing this back in the above relations one gets the relations
$$\alpha_1D_1+\alpha_2D_2+\alpha_3 D_3=0\quad {\rm and} \quad \beta_1D_1+\beta_2D_2+\beta_3 D_3=0 $$  thus forcing $\alpha_1=\alpha_2=\alpha_3=\beta_1=\beta_2=\beta_3=0$ as $\{D_1,D_2,D_3\}$ is $k$-linearly independent.

We have proved that $\dim_k(R_1J_d)\geq 6$.
We can do one better easily, to wit, force $\ell$ to avoid the prime $(x,y)$ as well.
We claim that for some $1\leq j\leq 3$, one has $yD_j\notin k\,xD_1+k\,xD_2+k\,xD_3+k\, \ell D_1+k\, \ell D_2+k\, \ell D_3$.
Supposing the contrary, 
in particular,  $yD_j\in (x,\ell)$ for $1\leq j\leq 3$. But since $\ell\notin (x,y)$ by our choice, then  $y\notin (x,\ell)$ and hence  $D_j\in (x,\ell)$ for every $j\in\{1,2,3\}$. 
It follows that $(x,\ell)$ is a minimal prime of $k[x,y,z]/(D_1,D_2,D_3)$, thus a minimal prime of $R/J$ as well, contradicting the choice of $\ell$.

This gives $\dim_k(R_1J_d)\geq 7$.

We still have one to go.
Thus, assume $\dim_k(R_1J_d)= 7$, i.e., the ideal $I$ has $2$ independent linear syzygies.
Then the syzygy matrix of $J$ would have the form

$$\varphi=\left(\begin{array}{c|ccccccccc}
\begin{array}{cc}
\ell_{1,1}&\ell_{1,2}\\
\ell_{2,1}&\ell_{2,2}\\
\ell_{3,1}&\ell_{3,2}
\end{array}& 
\begin{array}{cccc}
q_{1,1}&q_{1,2}&\ldots&q_{1,d-2}\\
 q_{2,1}&q_{2,2}&\ldots&q_{2,d-2}\\
 q_{3,1}&q_{3,2}&\ldots&q_{3,d-2}
\end{array}\\
\hline
\begin{array}{cc}
0&0\\
0&0\\
0&0\\
\vdots&\vdots\\
0&0
\end{array}&
\begin{array}{cccc}
h_{1,1}  & \kern-16pt h_{1,2}& \kern-10pt \ldots & h_{1,d-2}\\
h_{2,1} & \kern-16pt h_{2,2}& \kern-10pt \ldots & \kern-3pt h_{2,d-2}\\
h_{3,1} & \kern-16pt h_{3,2}& \kern-10pt \ldots & \kern-3pt h_{3,d-2}\\
\vdots & \kern-10pt \vdots&& \kern-3pt \vdots\\
h_{d-2,1} & \kern-6pt h_{d-2,2}& \kern-3pt\ldots & \kern-3pt h_{d-2,d-2}\\
\end{array}
\end{array}\right)$$

Now, note that for each $1\leq j\leq 3$, the generator $D_j$ can be taken to be the determinant of the submatrix of $\phi$ obtained by omitting the $j$th row.
This gives that the determinant of the right lower corner of size $(d-2)\times (d-2)$ is a factor of every $D_j$. 
But this is absurd as we are assuming $d\geq 5$ and the linear system $J_d$ has no fixed part.
Therefore,  $\dim_k(R_1J_d)\geq 8$ as contended.

\medskip

To complete the proof of the supplementary assertions of this item we explain the issue about the de Jonqui\`eres map.
From one side the structure of a plane de Jonqui\`eres map (\cite[Proposition 2.3]{HS}) tells us that its base ideal admits a unique linear syzygy up to scalars, i.e., $s=1$ in the notation of (a).

Conversely, assume that $s=1$.
We now focus our attention on the $(d-3)\times (d-1)$ submatrix  $\mathfrak{L}$  of the matrix $\varphi$ obtained by omitting the  first three rows.

We note that $\mathfrak{L}$ has the form $(\mathbf{0}\,|\, \mathfrak{L}')$ where
$\boldsymbol 0$ is the $(d-3)\times 1$ zero matrix and $\mathfrak{L}'$ is a $(d-3)\times (d-2)$ matrix with linear entries. The ideal of maximal minors of the latter matrix has codimension $\geq 2$ as it contains $I$.
It follows readily that the following complex
$$0\rar R\lar R^{d-2}\stackrel{\mathfrak{L}'}\lar R^{d-3}$$
is exact, where the entries of the leftmost map are the (signed, ordered) maximal minors of $\mathfrak{L}'$.

On the other hand, the syzygies $\mathfrak{K}$ of  $\mathfrak{L}=(\mathbf{0}\,|\, \mathfrak{L}')$ are $(1\,0\,\ldots\,0)^t$ and the syzygies of $\mathfrak{L}'$.
It follows that $\mathfrak{K}$ is a matrix with $2$ columns.
Therefore, by Proposition~\ref{Supplement-MAIN}, (iv) below below the syzygies of the ideal $I$ are generated by $2$ syzygies, i.e., $I$ is a perfect ideal.

Now apply Proposition~\ref{one-linear}.

\medskip

(b)  
The proof will consist in showing that $(J_{d+1})$ has maximal analytic spread ($=3$) and is linearly presented. Then the stated result follows from \cite[Theorem 3.2]{AHA}.

Since $(J_{d+1})$ is equi-homogeneous, its analytic spread is the dimension of the $k$-algebra $ k[J_{d+1}]$.
But $k[xD_1,yD_1,zD_1]$ is a $k$-subalgebra of $k[J_{d+1}]$, while $\{xD_1,yD_1,zD_1\}$ is obviously algebraically independent over $k$. Therefore, $\dim k[J_{d+1}]=3$.

To proceed, we assume for simplicity  that the Cremona map is not a de Jonqui\`eres map; the result is similar in the latter case and the changes are tiny.
We actually show that  $(J_{d+1})$ has a free linear resolution of the shape
$$0\rar R(-(d+3))^3 \lar R(-(d+2))^{d+7}\lar R(-(d+1))^{d+5}\lar (J_{d+1})\rar 0.$$
To see this, we argue as follows.
Note that $J=(J_d,B)$ and $(J_{d+1})=(R_1J_d, B)$ where $B$ is a $k$-vector basis of the direct supplement of ${\rm Im}(R_1\otimes_k J_d \lar J_{d+1})$.
By (a), it follows that $(J_{d+1})$ is generated by $9+d-4=d+5$ forms of degree $d+1$.
Furthermore, again by (a) any minimal syzygy of $J$ is the transpose of a vector $(q_1 \, q_2\, q_3\, \ell_4 \, \cdots \, \ell_{d-1})$,
where each $q_i$ is a $2$-form and each $\ell_j$ is a linear form.

Writing every $q_i$ as an element of $\fm:=(x,y,z)$, with linear coefficients, readily yields a linear syzygy of the above generators of $(J_{d+1})$ (the expression of $q_i$ as an element of $\fm$ is not unique, but it will turn out to be irrelevant which one is picked).
This gives a total of $d-2$ linear syzygies.

One has to add up the reduced Koszul relations of each one of the following pairs $\{xD_i, yD_i\}, \{xD_i, zD_i\}, \{yD_i, zD_i\}$, for $i=1,2,3$.
These contribute $9$ additional syzygies, giving a total of $d-2+9=d+7$ syzygies. Note that the latter syzygies are just the diagonal sum of three times the matrix of Koszul syzygies of $\fm$, in other words, they generate the syzygies of the $R$-module $\fm\oplus\fm\oplus\fm$.
 
On the other hand, one has the obvious free resolution 
\begin{equation}\label{3times-m}
0\rar R(-3)^3 \stackrel{\Psi}{\lar} R(-2)^9 \lar R(-1)^9\lar \fm\oplus\fm\oplus\fm\rar 0.
\end{equation}
Splicing, we get the following free complex
$$0\rar R(-(d+3))^3 \stackrel{\Psi(-d)}{\lar} R(-(d+2))^{9+d-2}\stackrel{\Phi}{\lar} R(-(d+1))^{d+5}\lar (J_{d+1})\rar 0.$$
Now use the Buchsbaum--Eisenbud acyclicity criterion: the ranks are easily seen to add up rightly; as to the heights of the required Fitting, one has  $\hht I_3(\Psi(-d))=\hht I_3(\Psi)=3$, as is obvious from (\ref{3times-m}), while for $\Phi$ one notes its shape as a $(d+5)\times (d+7)$ matrix:

$$\Phi=\left(\begin{array}{c|ccccccccc}
\begin{array}{ccc}
K(\fm)&&\\
& K(\fm)&\\
&& K(\fm)
\end{array}& 
\begin{array}{c}
\widetilde{\mathcal{L}}
\end{array}\\
\hline\\ 
\begin{array}{c}
\mathbf{0}
\end{array}&
\begin{array}{c}
\mathcal{L}
\end{array}
\end{array}\right),$$
where, to the left of the vertical line, $K(\fm)$ stands for the first Koszul matrix of $\fm$ and the empty slots are all zeros, whereas the matrix to the right of the vertical line is the $(d-2)\times (d+5)$ matrix of linear syzygies as obtained above from the original syzygies of $J$, the lower matrix $\mathcal{L}$ coinciding with the one described in (a) and the upper $\widetilde{\mathcal{L}}$ certain linear coordinates.

Note that the diagonal sum of $K(\fm)$ has only rank $6$, so it can only contribute this much to get a nonzero minor of size $(d+4)\times (d+4)$; forcefully, then one has the use a full maximal minor of the right side matrix.
It is clear that this patching allows for the construction of at least two $(d+4)$-minors with no proper common factor. This shows that $\hht I_{d+4}(\Phi)\geq 2$, as required for the acyclicity criterion.
\qed

\begin{Remark}\label{using-GGP}\rm
The result of item (b) also follows from \cite[Theorem 2.1 (ii)]{GGP}. To see this, one observes that in the notation there $\sigma=d+1$ as follows from Lemma~\ref{initial_degree}.
Moreover, the condition on how a straight line intersects the fat scheme translates here, via B\'ezout theorem, into the inequality $\mu_1+\mu_2\leq d$, a condition that is always satisfied by a proper homaloidal type as is the present case.

By that result the image of the birational map is a smooth surface and, moreover,
the surface is arithmetically Cohen--Macaulay (\cite[Proposition 2.4]{GGP}) (hence also projectively normal).
Since the present situation is more special it would be of some interest to have a direct algebraic proof of the Cohen--Macaulayness assertion.
\end{Remark}

\subsection{Further homological details}

As a sequel/tool to the results of  Theorem~\ref{MAIN}, we have the following  structurally finer points: 

\begin{Proposition}\label{Supplement-MAIN}
Keep the notation of {\rm Theorem~\ref{MAIN}} and, as introduced in the proof of the supplementary assertion of item {\rm (a)}, let $\mathfrak{L}$ denote the  $(d-4+s)\times (d-2+s)$ submatrix of  $\varphi$  omitting the  first three rows.
Then:
\begin{enumerate}
\item[{\rm (i)}] 
The $R$-module $J/I$ is linearly presented by $\mathfrak{L}$ if $s=0$ {\rm (}respectively, by its submatrix omitting the zero column if $s=1${\rm )}.
\item[{\rm (ii)}]  $J=I^{\rm sat} \Rightarrow \dim (J/I)=0 \Rightarrow \mu_i\leq 2,\;\forall i  \Rightarrow (d;\boldsymbol\mu)=(4;2^3,1^3)$ or $(d;\boldsymbol\mu)=(5;2^6) \Rightarrow J=I^{\rm sat}$, hence all these conditions are mutually equivalent.
\item[{\rm (iii)}] If $\dim (J/I)=1$ then $R/I$ is Cohen--Macaulay if and only if $J/I$ is Cohen--Macaulay.
\item[{\rm (iv)}]
Setting $\mathfrak{K}:=\ker(\mathfrak{L})$, the minimal syzygy matrix $\psi$ of $I$ is the submatrix of $\varphi\cdot \mathfrak{K}$ consisting of its first three rows. 
\end{enumerate}
\end{Proposition}
\demo
(i) By the part of  Theorem~\ref{MAIN},(a) which is independent of the assertion about the de Jonqui\`eres maps, we know that $s=1$ means that exactly one column of $\phi$ such that the entries on its first three rows are linear and the rest are zeros.
Since the generators of $I$ are minimal generators of $J$, the result is clear. 

\smallskip

(ii) 
Clearly, $\dim (I^{\rm sat}/I)=0$, so the first implication of the stated string is obvious.

To proceed, by (i) one has $\hht (I:J)=\hht I_{d-4+s}(\mathfrak{L})$ since $J/I$ and its zeroth Fitting ideal have the same radical.
Therefore, $\dim (J/I)=\dim (R/I:J)=3- \hht I_{d-4+s}(\mathfrak{L})$, so we are to show that $\hht I_{d-4+s}(\mathfrak{L})=3 \Rightarrow \mu_i\leq 2,\: \forall i$.
Assuming the contrary one has $\mu_1\geq 3$.
Then $d-1+s-\mu_1\leq d-1+s-3=d-4+s$, hence for even more reason $\hht I_{d-1+s-\mu_1}(\mathfrak{L})\geq 3$.
But letting $P_1$ denote the prime ideal for the first multiplicity, we see that $I_{d-1+s-\mu_1}(\mathfrak{L})\subset P_1$ since otherwise we could invert a minor in the localization $R_{P_1}$ and thus have $J_{P_1}={P_1}_{P_1}^{\mu_1}$ generated by $d-1+s-(d-1+s-\mu_1)=\mu_1$ elements -- this is absurd as ${P_1}_{P_1}^{\mu_1}$ is minimally generated by $\mu_1+1$ elements.

We move on to the next implication, namely that $\mu_i\leq 2,\: \forall i$ can only happen for those two particular homaloidal types.
This is an easy calculation using the equations of conditions.
Thus, letting $(d;2^{r_2}, 1^{r_1})$ be the assumed type, one has
$$
\left\{
\begin{array}{c}
r_1+2r_2=3d-3\\
r_1+4r_2=d^2-1
\end{array}
\right.
$$
Manipulation of these equations yields $r_1=d(6-d)-5$, from which (since $r_1\geq 0$) the only possibilities are $d=4, r_1=3$ or else, $d=5, r_1=0$; correspondingly, $r_2=3$ or else $r_2=6$.

The last implication is well-known -- for $d=4$ the base ideal $I$ is always saturated (cf. \cite[Theorem 1.5 (i)]{HS}), while $(5;2^6)$ is a classic (see, e.g., \cite[Theorem 2.14 (iii) Supplement]{HS} and the preliminaries of the proof, where it is argued that $(x,y,z)f\in I$, where $f$ is the unique minimal generator of degree $6$ of $J$).

\smallskip

(iii) This follows easily chasing depth along the following structural exact sequences
$$0\rar I \lar R \lar R/I\rar 0 \quad {\rm and}\quad 0\rar I \lar J\lar J/I\rar 0,$$
using that $\dim(J/I)=1$ by hypothesis and that depth$(J)=2$ and depth$(R)=3$.

\smallskip

(iv) Note that any syzygy of $\{D_1,D_2,D_3\}$ stacked on top of the $(d-4+s)\times 1$ zero matrix is a syzygy of $J$. Since $I=(D_1,D_2,D_3)$, one can write
$$\left(
\begin{array}{c}
\psi\\
\hline\\ [-10pt]
\underline{\mathbf{0}}
\end{array}
\right)=
\phi\cdot \mathfrak{C}
$$
for a suitable matrix $\mathfrak{C}$ (``content'' matrix), where the matrix on the left hand side is the matrix $\psi$ of syzygies of $I$ stacked on top of a zero matrix with $(d-4+s)$ rows and same number oc olumns as $\psi$ -- note that in the case where $s=1$ this content matrix will involve the unit column $(1\,0\,\ldots \,0)^t$.
Clearly, ${\rm Im}(\mathfrak{C})\subset \mathfrak{K}:=\ker(\mathfrak{L})$; in other words, every column of $\mathfrak{C}$ is a syzygy of $\mathfrak{L}$.

We claim the reverse inclusion as well.
To see this, let $\rho$ denote the submatrix of $\phi$ with the first three rows.
Then we wish to show that
$$(D_1\,D_2\,D_3)\cdot \rho\cdot \mathfrak{K}=0.$$
Recall that $D_j$ is the minor of $\phi$ omitting row $j$, for $j=1,2,3$.
Thus, if $\tilde{\phi}$ denotes the submatrix of $\phi$ omitting its last row and $\mathfrak{Q}$ denotes the matrix of cofactors of $\tilde{\phi}$, one has $(D_1\,D_2\,D_3)=L\cdot \mathfrak{Q}$, where $L$ denotes the last row of $\mathfrak{L}$.

Then one has

$$
(D_1\,D_2\,D_3)\cdot\rho\cdot \mathfrak{K}
= L\cdot\mathfrak{Q}\cdot\rho\cdot\mathfrak{K}
= \det(\tilde{\phi})\cdot L\cdot\mathfrak{K}=0
$$as was to be shown.
\qed

\section{Cremona maps with large highest multiplicity}

Some classical inequalities for Cremona map stem from the consideration of the three highest virtual multiplicities of the corresponding homaloidal type.
As far as we know, not much however has been obtained by stressing the role of the behavior of the map in the neighborhood of a single point with highest virtual multiplicity.

The goal of this part is to point out substantial algebraic information by focusing on when the highest virtual multiplicity $\mu$ roughly lies on the interval $d-1\geq \mu \geq\left\lceil \frac{d}{2} \right\rceil$, with a big question mark about the subinterval $d-5\geq\mu\geq \left\lceil \frac{d}{2} \right\rceil$.

Applying an arithmetic quadratic transformation (see \cite[Definition 5.2.3]{alberich}) to a homaloidal type belonging to this imprecise family will often yield a homaloidal type of the same shape (informally: ``telescopical''), while in a few other instances will afford a homaloidal type belonging to the inverse of the former.
Thereby, using repeatedly (or inductively) the Hudson test will show that these homaloidal types are proper.

\subsection{Classes with saturated base ideal}

The first two families below are known from the literature (see, e.g., \cite[Example 2.9]{Blanc}), but not so much the argument showing these are the only families with the chosen value of the highest virtual multiplicity, nor the question of the saturation of the base ideal.

\subsubsection{Highest virtual multiplicity $d-1$}

The case where the highest virtual multiplicity is also highest possible, namely, $\mu=d-1$, is well-known.
With this value of $\mu$, the remaining virtual multiplicities are forcefully all $1$'s (as it follows immediately from the equations of condition).
Thus, the corresponding homaloidal type is $(d;d-1,1^{2d-2})$.
Applying an arithmetic quadratic transformation based on the three higheste multiplicities immediately gives $(d-1;d-2,1^{2(d-1)-2})$, hence induction gives that the type is proper.
A Cremona map with this homaloidal type is the classical 
de Jonqui\`eres map.
The homological features of a de Jonqui\`eres are known (see, e.g., \cite[Theorem 2.7]{HS} and \cite[Corollary 2.6 ]{ST}).

\subsubsection{Highest virtual multiplicity $d-2$}

Assume $d\geq 4$ to avoid the de Jonqui\`eres map of degree $3$.

We prove:

\begin{Proposition}\label{highest_is_2less}  Let $T$ denote a proper homaloidal type of degree $d\geq 4$ and highest virtual multiplicity  $d-2$.
Then $T=(d;d-2,2^{d-2},1^{3})$.
Moreover,   given a  plane Cremona map $F$ of degree $d\geq 4$, with this homaloidal type and with proper base points, then $R/I$ is Cohen--Macaulay, where $I$ is the base ideal of $F$.
\end{Proposition}

\demo We first argue that any proper type $T$ with highest virtual multiplicity  $d-2$ is necessarily of the stated form.
 Applying an arithmetic quadratic transformation based on the highest virtual multiplicities gives the virtual multiplicity $d-(d-2)-\mu_2=2-\mu_2$. Since $T$ is proper, this must be nonnegative, i.e., $\mu_2\leq 2$. 
Thus, for a proper type necessarily $\mu_i\leq 2$ for every $i\geq 2$.
Therefore, we can set $T=(d;d-2, 2^s,1^{r-s})$.

It is straightforward to deduce the following relations from the equations of condition:
\begin{equation}
\label{two_less_eqs}
\left\{
\begin{array}{ccc}
 r+s &=& 2d-1\\
 r+3s&=&4d-5
\end{array}
\right.
\end{equation}
The solution is straightforward: $s=d-2, r=d+1$.

\medskip

To show that such a $T$ is proper we apply an arithmetic quadratic transformation based on the virtual multiplicities $d-2,2,2$ (recall that $d\geq 4$).
The result is the type $(d-2;d-4,2^{d-4},1^{3})$, hence the Hudson test works telescopically, yielding the main inductive step.
Therefore, it suffices to show that $(4;2; 2^{4-2}, 1^3)=(4;2^3,1^3)$ is a proper type to start the induction.
But this is obvious as an additional arithmetic quadratic transformation yields the homaloidal type $(2;1^3)$ of a quadratic map.

\medskip

Finally, we prove the last statement.
For it, let $p=p_1$ denote the point with virtual multiplicity $d-2$. 

Since $F$ is not de Jonqui\`eres, by Proposition~\ref{one-linear} the free resolution of the associated ideal $J$ of points with the virtual multiplicities of $F$ has the first of the two possible forms in  Theorem~\ref{MAIN}:

\begin{equation}\label{non-Jonq}
0\rar R(-(d+2))^{d-2}\stackrel{\varphi}\rar R(-(d+1))^{d-4}\oplus R(-d)^3\rar J\rar 0,
\end{equation}
where $I=(J_d)$ is the base ideal of $F$.

We claim that $I_1(\phi)\subset P$, where $P\subset R$ is the prime ideal of $p$.
Indeed, otherwise locally at $P$ one can invert an entry of $\phi$ so that by (\ref{non-Jonq}) the localization $J_P$ would be generated by at most $d-2$ elements. But this contradicts the fact that $J_P=P_P^{d-2}$ is minimally generated by $d-1$ elements because $P$ is generated by two linear forms (as we are assuming that $k$ is algebraically closed). 

This implies, in particular that $I_1(\mathfrak{L})\subset P$, where $\mathfrak{L}$ denotes the $(d-4)\times (d-2)$ linear submatrix of $\phi$, as in Proposition~\ref{Supplement-MAIN}.
Now, we may assume at the outset that $P=(x,y)$. Since the entries of $\mathfrak{L}$ are linear forms, then every such entry is actually a linear form in the ring $k[x,y]$.
Viewed as a module over $k[x,y]$ the kernel of $\mathfrak{L}$ is a free module, necessarily of rank $2$.
By flatness, this is also the nature of the kernel considered as a module over $R=k[x,y,z]$.

By Proposition~\ref{Supplement-MAIN} the minimal number of generators of the syzygy module of the base ideal $I$ of $F$ coincides with the minimal number of generators of the syzygy module of $\mathfrak{L}$.
This shows that the minimal number of generators of the syzygy module of $I$ is $2$, thus proving the statement.
\qed

\subsection{Classes with conjectured saturated base ideal}

\subsubsection{Highest virtual multiplicity $d-3$}

As in the previous case, we will face a strong restriction on the remaining virtual multiplicities.
However, differently from the previous case where one has a unique class of proper homaloidal types, here an enriched family will pop up,  subsuming subclasses distinguished according to the values of $d$ modulo $3$.
The typical proper homaloidal type that will emerge has as remaining virtual multiplicities $3,2,1$, with $d-3\geq 3$.

\begin{Proposition}\label{highest_is_3less}  Let $T$ denote a proper homaloidal type of degree $d\geq 6$ and highest virtual multiplicity  $d-3$.   
Then $T$ is one of the following, where $t\geq 2$:
\begin{enumerate}
\item[{\rm (0)}] $d\equiv 0\pmod 3$

Proper homaloidal types {\rm (}$d=3t${\rm )}:
\begin{enumerate}
\item[{\rm (0.1)}] $(d;d-3,3^{2t-2},2,1^4)$
\item[{\rm (0.2)}] $(d;d-3,3^{2t-3},2^4,1)$
\end{enumerate}
\item[{\rm (1)}] $d\equiv 1\pmod 3$

Proper homaloidal types {\rm (}$d=3t+1${\rm )}:
\begin{enumerate}
\item[{\rm (1.1)}] $(d;d-3,3^{2t-1},2^0,1^5)$
\item[{\rm (1.2)}] $(d;d-3,3^{2t-2},2^3,1^2)$
\end{enumerate}
\item[{\rm (2)}] $d\equiv 2\pmod 3$

Proper homaloidal types {\rm (}$d=3t+2${\rm )}:
\begin{enumerate}
\item[{\rm (2.1)}] $(d;d-3,3^{2t-1},2^2,1^3)$
\item[{\rm (2.2)}] $(d;d-3,3^{2t-2},2^5,1^0)$
\end{enumerate}
\end{enumerate}
{\rm (}{\sc Conjectured supplement{\rm )}}
Moreover,  given a  plane Cremona map $F$ of degree $d\geq 5$, with any of the stated homaloidal types and with proper base points, then $R/I$ is Cohen--Macaulay, where $I$ is the base ideal of $F$.
\end{Proposition} 
\demo 
We first argue that any proper type $T$ with highest virtual multiplicity  $d-3$ is necessarily of the form $T=(d;d-3, 3^{r_3}, 2^{r_2},1^{r_1})$.
Since $T$ is proper and $\mu_2\leq d-3$, then applying an arithmetic quadratic transformation based on the highest virtual multiplicities yields the virtual multiplicity $d-(d-3)-\mu_2=3-\mu_2$, which is nonegative only if $\mu_2\geq 3$. 
Thus, for a proper type necessarily $\mu_i\leq 3$ for every $i\geq 2$.
Therefore, we can set 
 $T=(d;d-3, 3^{r_3}, 2^{r_2},1^{r_1})$.

The equations of condition yield
$$\left\{
\begin{array}{ccc}
r_1+2r_2+3r_3+d-3 & =& 3d-3\\
r_1+4r_2+9r_3+(d-3)^2 &=&d^2-1
\end{array}
\right.
$$
An easy manipulation of these two relations yields the following ones
\begin{equation}
\label{three_less_eqs}
\left\{
\begin{array}{ccc}
 r_1+r_2 &=& 5\\
 r_2+3r_3&=&2d-5
\end{array}
\right.
\end{equation}
The solutions are straightforward by taking a nonnegative integer partition of $5$ and plugging into the second of these relations, thus producing a congruence modulo $3$, whose solutions are pretty immediate.

To see that all these homaloidal types are proper we apply the Hudson algorithm in terms of arithmetic quadratic transformations based on three highest multiplicities (cf. \cite[Corollary 5.2.21]{alberich}.
As above, in the case of highest multiplicity $d-1$ or $d-2$, the result of applying such an arithmetic transformation is telescopic.
Indeed, for $t\geq 2$  the arithmetic transformation is exactly the same in all cases, as the three highest multiplicities are $d-3,3,3$.
Moreover, since one goes down by $3$ in the degree every time an arithmetic quadratic transformation is applied, then the degree of the resulting type does not leave its class modulo $3$.
Explicitly, one obtains:

\medskip

\noindent Case $(0.1)$ $(d-3;d-6, 3^{2(t-2)},2,1^4)$

\noindent Case $(0.2)$ $(d-3;d-6,3^{2t-5},2^4,1)\;(t\geq 3)$

\noindent Case $(1.1)$ $(d-3;d-6, 3^{2t-3},2^0,1^5)$

\noindent Case $(1.2)$ $(d-3;d-6, 3^{2(t-2)},2^3,1^2)$

\noindent Case $(2.1)$ $(d-3;d-6, 3^{2t-3},2^2,1^3)$

\noindent Case $(2.2)$ $(d-3;d-6, 3^{2(t-2)},2^5,1^0)$

\medskip

To start the induction one needs the bottom type in each case, along with one additional arithmetic quadratic transformations $\mathcal{Q}$:

\medskip

\noindent Case $(0.1)$ $(6;3=6-3, 3^2,2,1^4)=(6;3^3,2,1^4)\; \stackrel{\mathcal{Q}}{\leadsto}\; (3;2,1^4)$: de Jonqui\`eres

\noindent Case $(0.2)$ $(6;3=6-3, 3,2^4,1)=(6;3^2,2^4,1)\; \stackrel{\mathcal{Q}}{\leadsto}\; (4;2^3,1^3)$: the non de Jonqui\`eres

\noindent Case $(1.1)$ $(7;4=7-3, 3^3,2^0,1^5)=(7;4,3^3,1^5)\; \stackrel{\mathcal{Q}}{\leadsto}\; (4;3,1^6)$: de Jonqui\`eres

\noindent Case $(1.2)$ $(7;4=7-3, 3^2,2^3,1^2)=(7;4,3^2,2^3,1^2)\; \stackrel{\mathcal{Q}}{\leadsto}\; (4;2^3,1^3)$: the non de Jonqui\`eres

\noindent Case $(2.1)$ $(8;5,3^3,2^2,1^3)\; \stackrel{\mathcal{Q}}{\leadsto}\; (5;3, 2^3,1^3)$: the saturated non de Jonqui\`eres

\noindent Case $(2.2)$ $(8;5,3^2,2^5)\; \stackrel{\mathcal{Q}}{\leadsto}\; (5; 2^6)$: the non-saturated.

\smallskip

This shows the properness of all types.

\begin{Remark}\rm
Note that the virtual multiplicity $2$ is present in all cases above with the exception of (1.1). Thus, in all these cases one has $\mu_1+2=d-3+2=d-1$. By \cite[Proposition 5.2]{Blanc} this means that when the highest multiplicity is $d-3$, except for case (1.1),  the Cremona map  belongs to the closure of the set of Cremona maps of degree $d+1$.
\end{Remark}

\subsubsection{Highest virtual multiplicity $d-4$}

We first note that any proper homaloidal type $T$ with highest virtual multiplicity  $d-4$ is necessarily of the form
$$T=(d;d-4, 4^{r_4}, 3^{r_3}, 2^{r_2}, 1^{r_1}).$$
Indeed, since $T$ is proper and $\mu_2\leq d-4$, then applying an arithmetic quadratic transformation based on the highest virtual multiplicities yields the virtual multiplicity $d-(d-4)-\mu_2=4-\mu_2$, which is nonegative only if $\mu_2\geq 4$. 
Thus, for a proper type necessarily $\mu_i\leq 4$ for every $i\geq 2$.

We will accordingly assume this form at the outset in the next proposition.  
In addition, the bulk of this case will have $d-4\geq 4$, i.e., $d\geq 8$, though we harmlessly include $d=7$ as a sort of degenerate case.

\begin{Proposition}\label{highest_is_4less}  Let $T=(d;d-4, 4^{r_4}, 3^{r_3}, 2^{r_2}, 1^{r_1})$ denote a proper homaloidal type of degree $d$ and highest virtual multiplicity  $d-4\geq 3$.  
One has:
\begin{enumerate}
\item[{\rm (a)}] $r_3$ and $r_4$ cannot vanish simultaneously.
\item[{\rm (b)}] $r_1$ and $r_2$ cannot vanish simultaneously.

\item[{\rm (c)}] If $r_4=0$ then $T$ is one of the following  types 
$$(7;3=7-4,3^3,2^3) =(7;3^4,2^3),\; (8;4=8-4,3^5,1^2),\; (10; 6=10-4, 3^7, 1^0),
$$
where the second is obtained from the last applying an arithmetic transformation on the three highest multiplicities.
\item[{\rm (d)}] 
If $r_3=0$ then  $r_4=\left \lfloor{(d-3)/2}\right \rfloor\geq 2$ and
$$T=\left\{
\begin{array}{ll}
(d;d-4, 4^{\left \lfloor{(d-3)/2}\right \rfloor}, 2^{3}, 1^{3}) & \mbox{if $d$ is even}\\ [6pt]
(d;d-4, 4^{(d-3)/2}, 1^{7}) &  \mbox{if $d$ is odd and $\not\equiv 3 \bmod 4$}
\end{array}
\right.
$$
\item[{\rm (e)}] Suppose that neither $r_3$ nor $r_4$ vanishes.
\begin{enumerate}
\item[{\rm (e.1)}] $r_1=0, r_2\neq 0$: in this case, $T=(d;d-4, 4^{(d-7)/2},3^3, 2^3)$, where  $d\geq 9$ is odd.
\item[{\rm (e.2)}] $r_2=0, r_1\neq 0$: in this case, $T=(d;d-4, 4^{(d-(r_3+3))/2}, 3^{r_3}, 1^{7-r_3})$, where $d$ is even and $r_3=1,3,5$.
\end{enumerate} 
\end{enumerate} 
{\rm (}{\sc Conjectured supplement{\rm )}}
Moreover,  given a  plane Cremona map $F$ of degree $d\geq 7$, with any of the stated homaloidal types and with proper base points, then $R/I$ is Cohen--Macaulay, where $I$ is the base ideal of $F$.
\end{Proposition}
\demo
Consider the equations of condition for this type:
$$\left\{
\begin{array}{ccc}
 r_1+2r_2+3r_3+4r_4 &=& 3d-3-(d-4)=2d+1\\ [6pt]
 r_1+4r_2+9r_3+16r_4&=&d^2-1-(d-4)^2=8d-17
\end{array}
\right.
$$
An easy manipulation of these relations allows us to express $r_1$ and $r_2$ in terms of $r_3$ and $r_4$ and $d$:
\begin{equation}\label{basic_relations}
\left\{
\begin{array}{ccc}
 r_1&=&3r_3+8r_4-4d+19\\ [6pt]
 r_2&=& 3(d-3-(r_3+2r_4))
\end{array}
\right.
\end{equation}
Since $r_1$ and $r_2$ are nonnegative, we have the following inequalities
\begin{equation}\label{inequalities}
3r_3+8r_4\geq 4d-19,\;r_3+2r_4\leq d-3.
\end{equation}

\medskip

(a) This one is clear: $r_3=r_4=0$ imply $r_1=19-4d<0$ for $d\geq 5$. 

\medskip

(b) Plugging $r_1=r_2=0$ into their expressions above yields two linear equations in $r_3,r_4$ in terms of $d$ which are readily solved, giving $2r_4=d-10, r_3=d-3-(d-10)=7$
Then $d$ must be even and $d\geq 12$ since $r_4>0$ for $r_1=r_2=0$ as will be proven in the next item.

Note the telescopical behavior:

{\Large
$$(d;d-4, 4^{\frac{d-10}{2}}, 3^7) \; \stackrel{\mathcal{Q}}{\leadsto} (d-4;d-8,4^{\frac{(d-4)-10}{2}}, 3^7)$$
}
while for $d=12$ one has  

{\large
$$ (12; 8,4,3^7) \; \stackrel{\mathcal{Q}}{\leadsto} (9;5, 3^6,1)\; \stackrel{\mathcal{Q}}{\leadsto} (7;3^5,1^2)\; \stackrel{\mathcal{Q}}{\leadsto} (5;3^2.1^6)\; \stackrel{\mathcal{Q}}{\leadsto} (3;-1,\ldots)
$$
}
Therefore, this homaloidal type is not proper.

\medskip

(c) If $r_4=0$ then $(4d-19)/3\leq r_3\leq d-3$.
Necessarily, $d-3\geq (4d-19)/3$   which means that $d\leq 10$.
Now, for $8\leq d\leq 10$ we then see that $r_3=d-3$.
Substituting for $r_3=d-3, r_4=0$  in the formulas of $r_1, r_2$ yields $r_1=10-d, r_2=0$.

The case where $d=10$ reduces to the case where $d=8$ by applying an arithmetic quadratic transformation, as one readily sees. Therefore, it suffices to check the latter. Applying again such a transformation yields the case $(0.1)$ of Proposition~\ref{highest_is_3less}, with $d=6$.

Now, as for $d=7$, the interval for $r_3$ allows for two possibilities: $r_3=3$ or $r_3=4$.
In the first case, one gets the homaloidal type $(7;3=7-4, 3^3, 2^3,1^0)=(7;3^4,2^3)$. Applying two successive arithmetic quadratic transformations yields the type $(3;2,1^4)$ -- a de Jonqui\`eres type.
As to the case $r_3=4$, we get the type  $(7;3=7-4, 3^4, 2^0, 1^3) =(7;3^4,1^3)$, then $(5;3^2, 1^6)$; the latter is not a proper type since for $d=5$ there are only three proper homaloidal types (de Jonqui\`eres, symmetric $(5; 2^6)$ and $(5;3,2^3,1^3)$).
 
Finally, the case $d=9$ reduces to the above non-proper type in degree $7$ by an arithmetic quadratic transformation, hence is not proper either.

\medskip

(d) With $r_3=0$, the inequalities (\ref{inequalities}) imply
$$\frac{4d-19}{8}\leq r_4\leq \frac{d-3}{2}.$$
But the positive difference between these two bounds is $19/8-3/2=7/8 <1$.
This yields the stated value of $r_4$ as the floor of half $d-3$.

On the other hand, the value of $r_2$ is now given by $r_2= 3(d-3-2r_4)$, from which follows that $r_2=3$ or $r_2=0$ according to whether $d$ is even or odd, respectively.

Finally, plugging each of these values into the expression  $r_1=8r_4-4d+19=\left \lfloor{(d-3)/2}\right \rfloor -4d+19$ readily gives the two respective values of $r_1$, as stated.

\medskip

To prove that both types are proper, one notes that since $d\geq 7$, then applying an arithmetic quadratic transformation to each of them yields a homaloidal type of the same shape by substituting $d-4$ for $d$.
Explicitly, since $r_4\geq 2$ then $d\geq 8$ and further the three highest virtual multiplicities are $d-4,4,4$. Thus, we get
\begin{eqnarray} \nonumber
 \lefteqn{(2d-(d-4)-4-4; d-(d-4)-4=0,d-(d-4)-4=0, d-4-4,  }  \kern2cm\\ \nonumber
&&  4^{\left \lfloor{(d-3)/2}\right \rfloor - 2}, \ldots) = (d-4;d-8,  4^{\left \lfloor{(d-4)-3)/2}\right \rfloor}, \ldots) 
\end{eqnarray} 
By recursion, according to the Hudson test, the two homaloidal types will be proper if and only if the resulting types of this operation for $7\leq d\leq 11$ are proper (this is because we jump $4$ at the time).

One has

\smallskip

\noindent $(11;7.4^4,1^7)\; \stackrel{\mathcal{Q}}{\leadsto} (7;4^2,3,1^7)\; \stackrel{\mathcal{Q}}{\leadsto} (3;-1,\ldots)$: non-proper

\noindent $(10;6,4^3,2^3,1^3)\; \stackrel{\mathcal{Q}}{\leadsto} (6;4,2^4,1^3)$: proper by Proposition~\ref{highest_is_2less} 

\noindent $(9; 5, 4^3, 1^7)\; \stackrel{\mathcal{Q}}{\leadsto} (5;4,1^8)$: de Jonqui\`eres

\noindent $(8;4^3, 2^3,1^3)\; \stackrel{\mathcal{Q}}{\leadsto} (4;2^3,1^3)$: the saturated non de Jonqui\`eres.

\smallskip

This show the statement of this item.

\medskip

(e.1) Let $r_1=0$. Then the first of the basic relations in (\ref{basic_relations}) yields $3r_3+6r_4=4d-19-2r_4$; upon substitution in the second relation, one obtains $r_2=2(r_4+5)-d$.
Since $r_2\neq 0$, $d<2(r_4+5)$. 
Back to the first relation gives $3r_3+8r_4<8(r_4+5)-19=8r_4+21$, from which we deduce that $r_3<7$.
But again the first of the relations in (\ref{basic_relations}) forces $4$ to divide $19+3r_3$. The only possibility in the interval $1\leq r_3\leq 6$ is $r_3=3$.
This in turn implies that $r_4=(d-7)/2$, so in particular $d$ is odd.
Also then $r_2=2((d-7)/2+5)-d=3$.

In conclusion, $T=(d;d-4, 4^{(d-7)/2}, 3^3,2^3)$.
To show it is proper, apply an arithmetic quadratic transformation to the three highest multiplicities; assuming $d\geq 11$ yields the same sort of type.
At the bottom one gets $(11;7;4^2,3^3.2^3)$ which in turns eventually transforms to $(5; 3,2^3,1^3)$, a well-known proper homaloidal type. 

\smallskip

(e.2) Let $r_2=0$. From the second relation $r_3+2r_4=d-3$ upon substitution on the first relation yields $r-1=10+2r_4-d$, and since $r_1\neq 0$, $d<10+2r_4$.
Back to the second relation gives once more $r_3<7$.
Also $r_4=(d-(r_3+3))/2$, so $r_3$ and $d$ must have opposite parities.
Further, $r_1=7-r_3$.

This gives the claimed form of $T$.
Changing by an arithmetic quadratic transformation is again telescopic, for $d\geq r_3+7$.
At the bottom, we find $T=(r_3+7; r_3+3, 4^2, 3^{r_3}, 1^{7-r_3})$. 
Letting $r_3$ run through $1$ to $6$ we easily test the obtained types for properness applying arithmetic quadratic transformations, thus finding that they are proper for $r_3$ odd and improper for $r_3$ even.

\begin{Remark}\rm
Again by \cite[Corollary 5.2]{Blanc}, except for the two cases in item (d), all such Cremona maps lie on the closure of Cremona maps of degree $d+1$.

Of course, this is pretty general: a proper homaloidal type whose  highest virtual multiplicity is, say,  $\mu_1=d-\ell$ divides into two classes of topological nature: the ones such that there exists a multiplicity $\mu_i=\ell-1$ correspond to Cremona maps of degree $d$ lying on the closure of the set of Cremona maps of degree $d+1$, while the rest is the union of the two classes where either $\mu_2=\ell$ and  $\mu_i\leq \ell-2$ for every $i\geq 3$ or else $\mu_i\leq \ell-2$ for  every $i\geq 2$, in which  neither case corresponds to Cremona maps of degree $d$ lying on the closure of the set of Cremona maps of degree $d+1$.
\end{Remark}

We close this part with a result showing that the line of argument used in the cases $\mu_1\leq d-2$ in order to prove saturation of the base ideal cannot be extended to the situation where $\mu_1\geq d-3$.
First, we isolate the following more technical result, using the notation introduced in the proof of Proposition~\ref{Supplement-MAIN} and of Proposition~\ref{highest_is_2less}:

\begin{Lemma}\label{minors_L} Consider a Cremona map with homaloidal type
$(d;\mu_1\geq\ldots\geq\mu_r)$ and proper base points. Letting $P_i\subset R$ denote the prime ideal of the proper base point $p_i$, one has
 $I_{d-1-\mu_i}(\mathfrak{L})\subset P_i,$ for every $\mu_i\geq 3$. In particular, the Fitting ideal $I_t(\mathfrak{L})=2$ has codimension $2$ in the range $d-1-\mu_i\leq t\leq d-4.$
\end{Lemma}
\demo If not, let there exist a $(d-1-\mu_i)$-minor of $\mathfrak{L}$ not belonging to some $P_i.$ 
This minor is invertible in $R_{P_i}$, hence the ideal $J_{P_i}\subset R_{P_i}$ is generated by $(d-1)-(d-1-\mu_i)=\mu_i$ elements, where $J\subset R$ denotes the corresponding ideal of fat points. But  $J_{P_i}={P_{i}}^{\mu_i}_{P_i}$ is generated by exactly $\mu_i+1$ elements, yielding a contradiction.
\qed

\smallskip

Now for the promised statement:

\begin{Proposition}\label{punchline} 
Consider a Cremona map with homaloidal type
$(d;\mu_1\geq\ldots\geq\mu_r)\neq (5;2^6)$ and proper base points. 
The following condition are equivalent:
\begin{enumerate}
\item[{\rm (a)}] $\mu_1\geq d-2$
\item[{\rm (b)}] $I_1(\mathfrak{L})$ has codimension $2.$
\end{enumerate} 
\end{Proposition}
\demo
The implication (a) $\Rightarrow$ (b) has already been explained in the proof of Proposition~\ref{highest_is_2less}.

For the reverse implication,  we may assume that $d\geq 5$ as otherwise (a) is automatically satisfied.
Note that $I_1(\mathfrak{L})$ is a prime ideal since it is generated by linear forms, and $I\subset I_1(\mathfrak{L}) $. Therefore, it is a minimal prime of $R/I$, say, $I_1(\mathfrak{L})=P_i$ for some $1\leq i\leq r.$ 
In particular, $I\subset I_{d-4}(\mathfrak{L})\subset P_i^{d-4}$, thus implying that  any form in $I_d$ goes through $p_i$ with multiplicity $\geq d-4$. Since $I_d$ defines a Cremona map, the general form in $I_d$ has multiplicity exactly $\mu_i$ on $p_i$. Therefore, $d-4\leq\mu_i$, hence also $d-4\leq \mu_1$.

Let us contradict the assumption that $\mu_1\leq d-3$.

Recall that one is assuming that the type is not $(5;2^6)$.
Then, by the classification of types in Proposition~\ref{highest_is_3less} and Proposition~\ref{highest_is_4less}, if  $d-4\leq \mu_1\leq d-3$ then at least two different virtual multiplicities are $\geq 3$.
If it happens that one of these is $\mu_i$, we choose the other one, say, $\mu_j$ ($j\neq i$).
Letting $P_j$ denote the correponding prime ideal, we have two distinct base primes $P_i$ and $P_j$.

By Lemma~\ref{minors_L}, the Fitting ideal $I_{d-1-\mu_j}(\mathfrak{L})$ has codimension $2$ and $P_j$ is necessarily a minimal prime thereof.

On the other hand, we contend that $I_{d-1-\mu_j}(\mathfrak{L})$ is $P_i$-primary, which gives the sought contradiction.
Actually, this is true of any $I_t(\mathfrak{L})$ such that $1\leq t\leq d-4$.

Indeed, first note that $I_t(\mathfrak{L})$  has codimension $2$ since it contains $I$.
Now, say, $P_i=(\ell_1,\ell_2)$, two independent linear forms.
Since the entries of $\mathfrak{L}$ are linear forms, the required operations to achieve these linear forms are all $k$-linear operations. This implies that a minor of arbitrary size is a polynomial on $\ell_1,\ell_2$ with coefficients in $k$.
Thus, $P_i$ is the extension to $R$ of the maximal ideal of the subring $k[\ell_1,\ell_2]\subset R$ and $I_t(\mathfrak{L})$ is the extension to $R$ of an $(\ell_1,\ell_2)$-primary ideal of $k[\ell_1,\ell_2]$.
Therefore, as an ideal of $R$, $I_t(\mathfrak{L})$ is a $P_i$-primary ideal, i.e., $P_i$ is its unique associated prime. 
\qed

\section{Cremona maps with small highest multiplicity}

\subsection{Preamble on non-saturated base ideals}

In the previous section we saw some evidence to the effect that non-saturated Cremona maps of degree $d$ and  highest virtual multiplicity $>d-5$ do not exist, with the single exception of $d=5$ and highest virtual multiplicity $2=d-3$.
Since being non-saturated means that the map has a non-perfect base ideal, it seems of some interest to find out if this question bears any impact to the geometric side of the picture.
In this part we essay a few steps in this direction.

We are actually compelled to pose the following question:

\begin{Question}\label{Q1}\rm
Let $(d;\mu_1,\ldots)$ denote a proper homaloidal type and let $F$ stand for a Cremona map with this type and general base points. If the base ideal of $F$ is not saturated then $\mu_1\leq \left \lfloor{d/2}\right \rfloor$.
\end{Question}
Note that, for a Cremona map,  $\mu_2\leq \left \lfloor{d/2}\right \rfloor$ always holds, so the relevant case is when $\mu_1 > \mu_2$.

In other words, the question is whether the inequality $\mu_1> \left \lfloor{d/2}\right \rfloor$ triggers saturation, which is a lot more ambitious than the conjectured result of the previous section.
The gist of the present section rests on considering two conditions, apparently apart from each other. The first, as suggested, is to the effect that the highest virtual multiplicity be less than about half the degree of the map; the second requires that all virtual multiplicities be even numbers.

The next subsection deals with the second of these conditions, setting it up on a larger environment.

\subsection{The role of the symbolic square}

In this part we will focus on a set of virtual multiplicities $\boldsymbol\mu$ satisfying certain relations which
 are naturally related to the equations of condition of a Cremona map.
Doubling these multiplicities will in fact satisfy the equations of condition.
This is a purely arithmetic procedure, but it has some impact to the corresponding ideals of fat points.

\begin{Definition}\rm
Let $\boldsymbol\mu=\{\mu_1,\ldots,\mu_r\}$ be a set of nonnegative integers  satisfying the following condition: there exists an integer $s\geq 2$ such that
\begin{equation}\label{eqs4birational}
\sum_{i=1}^r \mu_i=3(s-1) \quad {\rm and}\quad \sum_{i=1}^r \mu_i^2=s(s-1).
\end{equation}
Note that the first of the above relations implies that $s$ is uniquely determined.
We will say that $\boldsymbol\mu$ is a {\em  sub-homaloidal multiplicity set in degree $s$}.
\end{Definition}
The reason for this terminology will soon become clear.

\begin{Example}\rm
Let $\mu_i=\mu, \; \forall i$ -- so called {\em uniform} case. Then the only solutions of the above equations are:
\begin{equation}\label{uniform}
\left\{
\begin{array}{ccc}
\mu=1, & r=6, & s=3\\ [6pt]
\mu=3, & r=8, & s=9
\end{array}
\right.
\end{equation}
\end{Example}
The calculation is straightforward: the two relations $r\mu=3(s-1),\; r\mu^2=s(s-1)$ readily implies $s=3\mu$, which in turns yields the relation $r\mu=3(3\mu-1)$. The latter is equivalent to $\mu(9-r)=3$, from which the two stated solutions follow suit.

\begin{Example}\rm
Let $s\geq 3$ be an odd integer and let 
$$\mu_1=\cdots =\mu_4=(s-1)/2,\, \mu_5=\cdots=\mu_{s+3}=1
$$ 
(thus, $r=s-1+4=s+3$).

The verification is immediate: $4(s-1)/2+(s-1)=3(s-1)$, while $4(s-1)^2/4+s-1=s^2-s$.
\end{Example}

We will say that  a set of nonnegative integers $\boldsymbol\nu=\{\nu_1,\ldots,\nu_r\}$ is a {\em homaloidal multiplicity set} in degree $t$ if the corresponding type $(t;\boldsymbol\nu)$ is homaloidal, i.e., satisfies the classical equations of condition for a Cremona map.
The reason to emphasize this seemingly pedantic terminology is to be able to treat both multiplicity notions on the same foot and so also corresponding ideals with respect to the same set of points, in case the respective multiplicity sets have the same number.

Quite generally, the arithmetic as well as the geometric relation between the two are as follows:

\begin{Lemma}\label{doubling}
Let $\boldsymbol\mu=\{\mu_1,\ldots,\mu_r\}$ stand for a set of nonnegative integers and let $s\geq 2$ denote an integer.
The following conditions are equivalent:
\begin{enumerate}
\item[{\rm (i)}] $\boldsymbol\mu$ is a sub-homaloidal multiplicity set in degree $s$
\item[{\rm (ii)}] $2\boldsymbol\mu$ is a homaloidal multiplicity set in degree $2s-1$.
\end{enumerate}
In addition, fixing a set of points $\mathbf{p}$ with same cardinality as $\boldsymbol\mu$, then $I(\mathbf{p}, 2\boldsymbol\mu)=I(\mathbf{p}, \boldsymbol\mu)^{(2)}$ {\rm (}second symbolic power{\rm )}.
\end{Lemma}
\demo
The proof of the equivalence is immediate noting that $2(3s-3)=3(2s-1)-3$ and $4(s2-s)=(2s-1)^2-1$.
The supplementary assertion follows since points in $\pp^n$ are complete intersections and therefore the symbolic power of arbitrary order of the ideal of a point coincides with the respective ordinary power.
\qed

\smallskip

The following result shows a stronger tie between the two types above.

\begin{Proposition}\label{initial_degree_of_divided}
Let $\boldsymbol\mu=\{\mu_1\geq\cdots \geq\mu_r\}$ denote a sub-homaloidal multiplicity set in degree $s\geq 3$ and 
let $\mathbf{p}\subset \pp^2$ stand for a set of $r$ general points. Write $J:=I(\mathbf{p}, \boldsymbol\mu)$.
If the ``doubled'' homaloidal type $(2s-1;2\boldsymbol\mu)$ is proper then ${\rm indeg}(J)=s$.
\end{Proposition}
\demo
Since the points in $\mathbf{p}$ are general there is a Cremona map with homaloidal type $(2s-1;2\boldsymbol\mu)$ and these base points. Letting $\mathbb{J}:=I(\mathbf{p},2 \boldsymbol\mu)$, by Lemma~\ref{initial_degree} the initial degree of $\mathbb{J}$ is $2s-1$, in the above notation.

First it is clear that ${\rm indeg}(J)\geq s$, by using that  $\mathbb{J}=J^{(2)}$ from Lemma~\ref{doubling}.
Indeed, if ${\rm indeg}(J)\leq s-1$ then ${\rm indeg}(J^2)\leq 2s-2 <2s-1$ and hence also ${\rm indeg}(\mathbb{J})< 2s-1$ as $J^2\subset J^{(2)}$. But this would contradict the result of Theorem~\ref{MAIN}, which says that ${\rm indeg}(\mathbb{J})\geq 2s-1$.

For the reverse inequality, an easy calculation gives $e(R/J)= (1/2) (s+3)(s-1)$, and hence ${{s+2} \choose 2} - e(R/J) = (s+5)/2>0.$
On the other hand, for any $t$ one has
$$\dim_k J_t \geq \max\left\{ {{t+2} \choose {2}} - e(R/J) ,0\right\}.
$$
Therefore, $\dim_k J_s\geq (s+5)/2>0$.
In particular, ${\rm indeg}(J)\leq s$, hence equality must be the case, as stated.
\qed

\medskip

The main drive of this section aims at the case where the homaloidal multiplicity set afforded by a sub-homaloidal multiplicity set is proper.
The following examples show some erratic behavior of the second symbolic power of an ideal of fat points with a sub-homaloidal multiplicity set not affording a proper homaloidal multiplicity set.
The calculations were done with the help of the {\em Macaulay} system by taking points with random coordinates.

To avoid repetition, if $J$ is an ideal of fat points then $\mathbb{J}$ will denote its symbolic square.

\begin{Example}\rm
The multiplicity set $(4^2, 1^{10})$ is sub-homaloidal in degree $7$, as one easily checks.  Here the doubled homaloidal type $(13;8^2,2^{10})$ is not proper as $13<8+8$. Set $J=I(\mathbf{p}; 4^2, 1^{10})$ for general points $\mathbf{p}$.
Then:
\begin{itemize}
\item $J$ is minimally generated in degrees $7$ and $8$ and  the subideal $(J_7)$ has codimension one
\item $\mathbb{J}$ is generated in degrees $13,14,15,16$ 
\item The ideal $(\mathbb{J}_{13})$ has codimension one.
\end{itemize} 
\end{Example}

\begin{Example}\rm
The multiplicity set $(4, 2^6, 1^2)$ is sub-homaloidal in degree $7$, as one easily checks.  As it turns, applying iteratively the Hudson test shows that the doubled homaloidal type $(13;8,4^6,2^2)$ is not proper. Setting $J$ and $\mathbb{J}$ as in the previous example, one gets:
\begin{itemize}
\item $J$ is minimally generated in the single degree $7$
\item $\mathbb{J}$ is generated in degrees $13,14$ and its free resolution of $\mathbb{J}$ has same Betti numbers and twists as (\ref{res_DJ}); in particular, $\mathbb{J}_{13}$ is a net
\item However, the ideal $(\mathbb{J}_{13})$ still has codimension one.
\item Canceling the $\gcd$ of the forms in $\mathbb{J}_{13}$ (the fixed part of the net) results in  a non-homaloidal net generating a saturated ideal.
\end{itemize} 
\end{Example}

\begin{Example}\rm 
The multiplicity set $(4, 3^6, 1^2)$ is sub-homaloidal in degree $9$ and once again, applying iteratively the Hudson test yields that the doubled homaloidal type $(17;8,6^6,2^2)$ is not proper. 
With same notation, one has:
\begin{itemize}
\item $J$ is minimally generated in the single degree $9$
\item $\mathbb{J}$ is generated in degrees $17,18$ and its free resolution of $\mathbb{J}$ has same Betti numbers and twists as (\ref{res_DJ}); in particular, $\mathbb{J}_{17}$ is a net
\item However, the ideal $(\mathbb{J}_{17})$ still has codimension one.
\item Canceling the $\gcd$ of the forms in $\mathbb{J}_{17}$ (the fixed part of the net) results in  a non-homaloidal net generating a non-saturated ideal.
\item However, the latter non-saturated ideal fails for the template of Proposition~\ref{resolving_non_saturated} in that it admits minimal syzygies in two different degrees.
\end{itemize} 

\end{Example}

\subsection{Non-saturated base ideal: the underlying geometry}

An additional remarkable feature of a Cremona map $\mathfrak{F}:\pp^2\dasharrow {\pp^2}^*$ is its {\em characteristic matrix}, which encapsulates the main geometric properties of the map.
On a pinhead description, the rows of the matrix, excluding the respective first entries, give the homaloidal type of the map and the types defining the principal curves corresponding to base points.
This part of the paper will indirectly evoke features of the characteristic matrix -- such as the notion of a principal curve -- but actually no essential use will be made of its theory, except for a basic criterion to be found in \cite{alberich} --  our standing reference for an exhaustive account of this theory.

From now on we will consider the situation of a sub-homaloidal multiplicity set in degree $s$ for which the highest multiplicity of the associated doubled homaloidal type is $s-1=\left \lfloor{(2s-1)/2}\right \rfloor$, as suggested at the beginning of the section. 

We will use the terminology and notation of the previous subsection.
Thus, $s\geq 3$ denotes an odd number.

\begin{Lemma}\label{proper_exceptional}
Let $(((s-1)/2)^3,\mu_4,\ldots,\mu_r)$ denote a sub-homaloidal multiplicity set in degree $s\geq 3$. If the associated ``doubled'' homaloidal set $((s-1)^3,2\mu_4,\ldots,2\mu_r)$ gives a proper homaloidal type in degree $2s-1$ then $(s-1; ((s-1)/2)^2,(s-3)/2,\mu_4,\ldots,\mu_r)$ is a proper exceptional type.
\end{Lemma}
\demo
Let us denote
$$T:=(2s-1;(s-1)^3,2\mu_4,\ldots,2\mu_r)$$
and
 $$S:=(s-1; ((s-1)/2)^2,(s-3)/2,\mu_4,\ldots,\mu_r),$$
 where one may harmlessly assume that $\mu_4\geq \cdots \geq\mu_r>0$.
 
 We first show that $S$ is an exceptional type.
 Namely, since $T$ is homaloidal then one has:
 $$\left\{\begin{array}{cccc}
 3(s-1)+2\sum_{i\geq 4}\mu_i &=& 6(s-1)\\ [4pt]
 3(s-1)^2+4\sum_{i\geq 4}\mu_i^2 &=& 4s(s-1),
 \end{array}\right.$$
 from which we derive $\sum_{i\geq 4}\mu_i =3(s-1)/2$ and $\sum_{i\geq 4}\mu_i^2=(s^2+2s-3)/4$.
 
 It follows that
$$\left\{\begin{array}{cccc}
 2(s-1)/2+ (s-3)/2+\sum_{i\geq 4}\mu_i &=& 3(s-1)-1\\ [4pt]
 2((s-1)/2)^2+((s-3)/2)^2+\sum_{i\geq 4}\mu_i^2&=& (s-1)^2-1,
 \end{array}\right.$$
 as required. 
 
Next, we look at properness.
For $i\geq 0$, define $T_i$ and $S_i$ iteratively by $T_0=T,S_0=S$ and for $i\geq 1$, $T_i$ (respectively, $S_i$) is the type obtained from $T_{i-1}$ (respectively, from $S_{i-1}$) by applying an arithmetic quadratic transformation based on three major multiplicities.

Since we are assuming that $T$ is a proper homaloidal type, then $T_j=(1;\underbrace{0,\ldots,0)}_r)$ for a certain $j\geq 1$.
Let $T_i$ be any transformed type, with $i<j$.

{\sc Claim:} $T_i$ is of the form $(2t;2\eta_1,\ldots,2\eta_\sigma,1^3,0^{r-\sigma-3}),$ for suitable $t\geq 1$, $\eta_1\geq\cdots\geq\eta_{\sigma}\geq 1$ and some $0\leq \sigma \leq r-3 .$

One proceeds by induction on $i\geq 1$;  applying an  arithmetic quadratic transformation based on the first three multiplicities, it is easily seen that $T_1=(s+1;2\mu_4,\ldots,2\mu_r,1^3).$
Note that now that the degree is even and the three highest multiplicities have been  shifted to the end;  after that, in every given stage at least the smallest of the three highest multiplicities will be shifted to give place to a higher {\em even} multiplicity.

Assume inductively that $T_{i-1}$ has the form $(2t;2\eta_1,\ldots,2\eta_\sigma,1^3,0^{r-\sigma-3})$ 
and let us show that necessarily $\sigma\geq 3$.
Indeed,  first is obvious that $\sigma\neq 0$ as otherwise  $T_{i-1}=(2;1^3)$, hence $T_i=(1;0,\ldots, 0)$, contradicting our assumption.
If $\sigma=1$ then the first equation of condition gives $2\eta_1+3=3.2t-3=6t-3$. In addition, one knows that the sum of the highest multilicities does not exceed the degree, hence $2\eta_1+1\leq 2t$. Together this yields $t\leq 1$, thus implying $\eta_1=0$; this contradicts the assumption that $\eta_1>0$.

The case where $\sigma=2$ is similar: $t\geq \eta_1+\eta_2=3t-3$ still yields $t=1$, forcing $\eta+1+\eta_2\leq 1$, i.e., $\eta_2=0$ -- a contradiction.

Now, since $\sigma\geq 3$, one sees that the transformed type is
$$(2(2t-\eta_1-\eta_2-\eta_3); 2(t-\eta_2-\eta_3) ,  2(t-\eta_1-\eta_3),  2(t-\eta_1-\eta_2), 2\eta_4,\ldots,2\eta_{\sigma}, 1^3, 0\ldots,0),$$
where the first three multiplicities are in decreasing order, but need to be reordered as compared to the subsequent multiplicities $2\eta_4\geq\cdots \geq2\eta_{\sigma}$.
In fact, some of the first three written above may turn out to vanish, thus increasing the number of zeros and decreasing $\sigma$.
However, since for $i\geq 1$ both the degree and the multiplicties are even, the number of $1$'s will not change.
The number of zeros is a consequence of the fact that the total number of multiplicities does not change under an arithmetic quadratic transformation.

This proves the claim.

\smallskip

Next, according to \cite[Proposition 5.5.5]{alberich}, an exceptional type is proper if and only it is ultimately transformable into a type $(0;-1, \underbrace{0,\ldots,0}_{r-1})$ after finitely many  arithmetic quadratic transformations.

{\sc Claim:} If $T_i=(2t;2\eta_1,\ldots,2\eta_\sigma,1^3,0^{r-\sigma-3})$ then $S_i=(t;\eta_1,\ldots, \eta_\sigma,1^2,0^{r-\sigma-2})$.

The proof is again by induction on $i$ and similar to the previous one, so we only mark the main points.
Thus, one easily sees that  $S_1=((s+1)/2;\mu_4,\ldots,\mu_r,1^2, 0).$
Supposing that $S_{i-1}=(t;\eta_1,\ldots, \eta_\sigma,1^2,0^{r-\sigma-2})$ the next stage will have the form
$$(2t-\eta_1-\eta_2-\eta_3; t-\eta_2-\eta_3, t-\eta_1-\eta_3,  t-\eta_1-\eta_2, \eta_4,\ldots,\eta_{\sigma}, 1^2,0\ldots,0).$$
This time around there may appear more $1$'s, but the number of zeros is always the number of zeros in $T_{i+1}$ plus one.
So much for the claim.

In particular, since in the one before the last stage of transforming $T$  one gets the type $(2;1^3,0^{r-3})$, then the corresponding transformed type of $S$ is $(1;1^2,0^{r-2})$ -- the latter clearly transforms to $(0;-1,0^{r-1})$, as was to be shown.
\qed

The main result of this part is as follows.

\begin{Theorem}\label{Main-non-saturated}
Let $(((s-1)/2)^3,\mu_4,\ldots,\mu_r)$ denote a sub-homaloidal multiplicity set in degree $s\geq 3$. Assume that the associated ``doubled'' homaloidal set $((s-1)^3,2\mu_4,\ldots,2\mu_r)$ gives a proper homaloidal type in degree $2s-1$. Let $I\subset R=k[x,y,z]$ denote the base ideal of a Cremona map $\mathfrak{F}$ on $r$ general points having this homaloidal type.
Then there are three principal curves of  $\mathfrak{F}$ of degree $s-1$ and three independent linear forms $\ell_1,\ell_2,\ell_3\in R$
such that $I=(\ell_3f_1f_2, \ell_2f_1f_3,\ell_1f_2f_3)$.
\end{Theorem}
\demo
Let $\{P_1,\ldots,P_r\}\subset R$ denote the defining prime ideals of $r$ general points over which $\mathfrak{F}$ is a Cremona map with homaloidal type  $T:=(2s-1;((s-1)^3,2\mu_4,\ldots,2\mu_r)$.
We may assume that the first three points are $(0:0:1), (0:1:0), (1:0:0)$.
From Theorem~\ref{MAIN} we know that the base ideal $I\subset R$ of  $\mathfrak{F}$ is generated by the homogeneous piece of degree $2s-1$ of the ideal of fat points
$$\mathbb{J}:=J^{(2)}=P_1^{s-1}\cap P_2^{s-1}\cap P_3^{s-1}\cap P^{2\mu_4}\cap\cdots\cap P^{2\mu_r},$$
where $J:=P_1^{(s-1)/2}\cap P_2^{(s-1)/2}\cap P_3^{(s-1)/2}\cap P^{\mu_4}\cap\cdots\cap P^{\mu_r}$.

Consider the following additional ideals of fat points:

\begin{equation}\label{jota1}
J_1=P_1^{(s-3)/2}\cap P_2^{(s-1)/2}\cap P_3^{(s-1)/2}\cap P_4^{\mu_4}\ldots\cap P_r^{\mu_r},
\end{equation}

\begin{equation}\label{jota2}
J_2=P_1^{(s-1)/2}\cap P_2^{(s-3)/2}\cap P_3^{(s-1)/2}\cap P_4^{\mu_4}\ldots\cap P_r^{\mu_r}
\end{equation}

\begin{equation}\label{jota3}
J_3=P_1^{(s-1)/2}\cap P_2^{(s-1)/2}\cap P_3^{(s-3)/2}\cap P_4^{\mu_4}\ldots\cap P_r^{\mu_r}.
\end{equation}

By Lemma~\ref{proper_exceptional}, the respective types of $J_1,J_2,J_3$ are proper exceptional types in degree $s-1$.
By \cite[Proposition 5.5.13]{alberich}, for each $i=1,2,3$, there exists a nonzero irreducible homogeneous polynomial $f_i\in J_i$ of degree $s-1$ defining a principal curve of  $\mathfrak{F}$.
Since we are assuming that $P_1=(x,y), P_2=(x,z), P_3=(y,z)$, the forms
 $xf_1f_2, zf_1f_3, yf_2f_3$ belong to $\mathbb{J}$.
But they are of degree $2s-1$, hence one has
 $(xf_1f_2, zf_1f_3, yf_2f_3)\subset I$.
On the other hand, no two $f_i,f_j (i\neq j)$ have a proper common factor as otherwise $f_i\in J$, while the latter has initial degree $s$ by Proposition~\ref{initial_degree_of_divided}.
Therefore, e.g., $\{f_1(xf_2+yf_3), zf_2f_3\}$ have no proper common factors, hence the ideal $(xf_1f_2, zf_1f_3, yf_2f_3)$ has codimension $2$.
By the same token, $\{xf_1f_2, zf_1f_3, yf_2f_3\}$ are $k$-linearly independent.
The equality now follows.
\qed

\begin{Proposition}\label{resolving_non_saturated}
Notation as in {\rm Theorem~\ref{Main-non-saturated}}. Then 
\begin{enumerate}
\item[{\rm (a)}]
$I$ is a non-saturated ideal with minimal graded resolution 
$$0\rar R(-3s)\stackrel{\psi}{\lar}R^3(-(3s-1))\stackrel{\varphi}{\lar} R^3(-(2s-1))\rar I\rar 0,$$
where 
$$\varphi=\left(\begin{array}{ccc}
\ell_2f_3&-\ell_3f_3&0\\
-\ell_1f_1&0&\ell_3f_1\\
0&\ell_1f_2&-\ell_2f_2
\end{array}\right)\quad \mbox{and}\quad \psi=\left(\begin{array}{c}\ell_3\\\ell_2\\\ell_1\end{array}\right)$$
\item[{\rm (b)}]
$I^{\rm sat}=(I,f_1f_2f_3).$
\end{enumerate}
\end{Proposition}
\demo
(a) The maps are easily seen to give a complex.
A simple application of the Buchsbaum--Eisenbud acyclicity criterion shows the contention.

(b) 
Since $I$ is the base ideal of a Cremona map, one has $I^{\rm sat}_{2s-1}=I_{2s-1}$ and $I^{\rm sat}_d=0$ for $d<2s-1$.
By part (a) and \cite[Proposition 1.3 (ii)]{HS}, and as $R/I^{\rm sat}$ is Cohen--Macaulay, it follows that $I^{\rm sat}$ is minimally generated by $4$ forms, three of which can be taken to be the minimal generators of $I$ as in (a) and the fourth one some form of degree $3(s-1)$.

Consider the following matrix
$$\left(\begin{array}{ccc}
f_1 & 0 & 0\\
0 & f_2 & 0\\
0 & 0 & f_3 \\
\ell_3 & \ell_2 & \ell_1
\end{array}\right)
$$
Its $3$-minors form a minimal set of generators of the ideal $(I,f_1f_2f_3)$, since for any $i\neq j$ one has $\gcd(f_i,f_j)=1$.
This ideal is contained in $I^{\rm sat}$ because $(x,y,z)f_1f_2f_3=(\ell_1,\ell_2,\ell_3)f_1f_2f_3\in I$.
Since $\deg f_1f_2f_3=3(s-1)$, we are done.
\qed

\noindent {\bf Authors' addresses:}

\medskip

\noindent {\sc Zaqueu Ramos},  Departamento de Matem\'atica, CCET\\ Universidade Federal de Sergipe\\
49100-000 S\~ao Cristov\~ao, Sergipe, Brazil\\
{\em e-mail}: zaqueu@gmail.com\\

\noindent {\sc Aron Simis},  Departamento de Matem\'atica, CCEN\\ Universidade Federal
de Pernambuco\\ 50740-560 Recife, PE, Brazil\\
{\em and}\\
 Departamento de Matem\'atica, CCEN\\
  Universidade Federal
da Para\ii ba\\
58059-900 Jo\~ao Pessoa, PB, Brazil.\\
{\em e-mail}:  aron@dmat.ufpe.br

\end{document}